\documentclass{IEEEtran}
\usepackage{cite}
\usepackage{amsmath,amssymb,amsfonts}
\usepackage{algorithmic}
\usepackage{graphicx}
\usepackage{textcomp}

\usepackage{epsfig} 
\usepackage{comment}
\usepackage{bm}
\usepackage{color}
\usepackage{algorithm}
\newtheorem{problem}{Problem}
\newtheorem{proposition}{Proposition}

\newtheorem{lemma}{Lemma}
\newtheorem{theorem}{Theorem}

\newtheorem{assumption}{Assumption}
\newtheorem{rmk}{Remark}

\newcommand{\argmin}{\mathop{\rm arg~min~}\limits}
\newcommand{\esssup}{\mathop{\rm ess~sup}\limits}

\usepackage[subrefformat=parens]{subcaption}
\def\BibTeX{{\rm B\kern-.05em{\sc i\kern-.025em b}\kern-.08em
    T\kern-.1667em\lower.7ex\hbox{E}\kern-.125emX}}
\begin{document}
\title{Non-convex optimization problems for \\ maximum hands-off control}
\author{Takuya Ikeda
\thanks{This study was supported in part by 
JSPS KAKENHI under Grant Number JP21K14188.}
\thanks{T. Ikeda is with Faculty of Environmental Engineering, 
The University of Kitakyushu, Fukuoka, 808-0135, Japan
(e-mail: t-ikeda@kitakyu-u.ac.jp). }}

\maketitle

\begin{abstract}
The maximum hands-off control is the optimal solution to the $L^0$ optimal control problem.
It has the minimum support length among all feasible control inputs. 
To avoid computational difficulties arising from its combinatorial nature, the convex approximation method that replaces the $L^0$ norm by the $L^1$ norm in the cost function has been employed on standard.
However, this approximation method does not necessarily obtain the maximum hands-off control.
In response to this limitation, this paper newly introduces a non-convex approximation method and formulates a class of non-convex optimal control problems that are always equivalent to the maximum hands-off control problem.
Based on the results, this paper describes the computation method that quotes algorithms designed for the difference of convex functions optimization.
Finally, this paper confirms the effectiveness of the non-convex approximation method with a numerical example.
\end{abstract}


\section{Introduction}

\IEEEPARstart{T}he theory of sparse representation is finding important applications in various fields of science and engineering~\cite{ZhaXuYan15}. 
The methodology tries to find a linear combination of a small number of basis vectors in a suitable space that better approximates an object vector.
To obtain such representation, various optimization methods have been proposed.

The natural penalty function for promoting sparsity is the $\ell^0$ norm that counts the number of nonzero components in a vector. 
The minimization of the $\ell^0$ norm is a combinatorial optimization problem known to be NP-hard~\cite{Nat95}.
To reduce the computational burden, a convex approximation method using the $\ell^1$ norm instead of the $\ell^0$ norm has been widely employed~\cite{CanPla09}. 
However, it has been reported that this method based on the $\ell^1$ minimization tends to cause bias to the estimate value due to the convexity of the penalty function.
To cope with this deficiency, many non-convex penalty functions have been proposed as seen in 
the smoothly clipped absolute deviation (SCAD)~\cite{FanLi01}, 
the minimax concave penalty (MCP)~\cite{Zha10}, 
the $\ell^p$ penalty with $0<p<1$~\cite{FraFri93}, 
and the log-sum penalty (LSP)~\cite{CanWakBoy08}.
Although for the non-convex penalty, it is very difficult to minimize in general, 
they are reportedly more effective theoretically and experimentally than the $\ell^1$ penalty in approximating the original $\ell^0$ penalty~\cite{Cha07,ChaSta08,TrzMan08,TraWeb19,YinMirPal20,WooCha16,YinLouHe15,SouBlaAub15}. 
For example, 
the studies~\cite{Cha07,ChaSta08,TrzMan08} on the $\ell^p$ optimization show
less-restrictive isometry conditions compared with the $\ell^1$ optimization,
and the possibility of precise reconstruction of sparse signals with fewer sample measurements.
The study~\cite{TraWeb19} derives the recovery guarantee of sparse signals by means of non-convex optimization, which is less restrictive than the standard null space property for the $\ell^1$ optimization.
%
The study~\cite{YinMirPal20} analyzes sparse graph learning in the Laplacian constrained Gaussian graphical model,
and clarifies an unexpected behavior of the $\ell^1$ regularization, i.e., learning a dense graph for a large regularization parameter.
Consequently, the study proposes a non-convex penalized maximum likelihood method that enables a theoretical guarantee against the statistical error.

The notion of sparsity has also been applied to dynamical control systems as seen, for example, in
sensor/actuator selection~\cite{ManKutBru21}, resource-aware control~\cite{GomHee15}, and state estimation~\cite{CarRusHes15}.
As the study most relevant to this paper, the optimal control to minimize the $L^0$ norm of control input in continuous-time systems is proposed in~\cite{NagQueNes16}.
Here, $L^0$ norm is a penalty function for measuring the length of support of the function.
This novel control is called the \emph{maximum hands-off control}.
The control characteristically allows actuators to be at a stop for a long time, thus contributing to the significant reductions in fuel consumption, power usage, and communication burden~\cite{Nag2020}.
Utilizing the idea of $\ell^1$ approximation in the sparse representation, 
the study~\cite{NagQueNes16} reveals the equivalence between the original $L^0$ optimal control problem and its convex approximation, $L^1$ optimal control problem, under an assumption called \emph{normality}. 
This result has been extended to general linear systems in~\cite{ITKKTAC18}, and some relevant studies have been published~\cite{ChaNagQueRao16,NagOstQue16,ItoIkeKas21}. 
%
%
On the other hand, if the normality assumption does not hold, the $L^1$ approximation method cannot necessarily obtain a sparse control as shown in~\cite{AthFal}.
Also, the usefulness of non-convex penalty functions for promoting sparsity has been reported as mentioned above.
Such a background motivates the analysis of non-convex optimal control problems in which less restrictive equivalence can be realized.

The main aimed contribution of this paper is to give a class of non-convex optimal control problems equivalent to the maximum hands-off control problem for continuous-time linear systems.
Unlike the $L^1$ approximation method, it is proved that this equivalence always holds.
For example, this class includes the optimal control problems in which the cost function is defined by the SCAD, the MCP, the LSP, and the $L^p$ norm with $0<p<1$.
In addition, it is also proved that some problems of this class can be represented as a difference of convex functions (DC) optimization problem.
As such, quoting the best-known DC algorithm, this paper describes a numerical algorithm for the maximum hands-off control.
Furthermore, it is confirmed that the maximum hands-off control can be obtained by solving the non-convex optimal control problems with an example where the normality assumption does not hold.
 
The remainder of this paper is organized as follows:
Section~\ref{sec:formulation} defines the maximum hands-off control problem and the non-convex optimal control problems, and introduces some examples of the non-convex optimal control problems;
Section~\ref{sec:analysis_comp} proves the equivalence as the main result of this study, and provides a computational algorithm;
Section~\ref{sec:example} confirms the results with a numerical example;
and Section~\ref{sec:conclusion} concludes the study.

\subsection*{Notations}\label{sec:math}

The set of all positive integers is denoted by $\mathbb{N}$,
and the set of all real numbers by $\mathbb{R}$.
For any $m\in\mathbb{N}$ and $\Omega\subset\mathbb{R}$,
$a=[a_1, a_2, \dots, a_m]^{\top}\in \Omega^{m}$ means that $a_i\in \Omega$ holds for all $i$.
The {\em $\ell^p$ norm} of $a\in\mathbb{R}^m$ is defined by
\begin{align*}
	&\|a\|_{\ell^0} \triangleq \#\{i\in\{1, 2, \dots, m\}: a_i\neq0\},\\
	&\|a\|_{\ell^p} \triangleq \left( \sum_{i=1}^{m} |a_i|^p\right)^{\frac{1}{p}}  \quad p\in(0, \infty),
\end{align*}
where $\#$ denotes the number of elements of the set.
For a matrix (or a vector) $M$, $M^{\top}$ denotes the transpose of $M$.
For $b, c\in\mathbb{R}$, $\min\{b, c\}$ represents the smaller value between $b$ and $c$.
For sets $S_1$, $S_2\subset \mathbb{R}$, 
$S_1 \backslash S_2\triangleq\{a \in S_1: a\not\in S_2\}$.
The Lebesgue measure on ${\mathbb{R}}$ is denoted by $\mu$. 
For a continuous-time signal $u(t)=[u_1(t), u_2(t), \dots, u_m(t)]^{\top}\in{\mathbb{R}}^m$ over a time interval $[0, T]$, its {\em $L^p$ norm} is defined by
\begin{align*}
	&\|u\|_{L^0} \triangleq \sum_{j=1}^{m}\mu(\{t\in[0, T]: u_j(t)\neq0\}),\\
	&\|u\|_{L^p} \triangleq \left\{\sum_{j=1}^{m}\int_{0}^{T} |u_j(t)|^p dt\right\}^{\frac{1}{p}} \quad p\in(0,\infty),\\
	&\|u\|_{L^\infty} \triangleq \max_{1 \leq j \leq m} \esssup_{0\leq t\leq T} |u_j(t)|,
\end{align*}
where $T>0$.
The set of all functions that have the finite $L^p$ norm on a measurable set $E\subset\mathbb{R}$ is denoted by $L^p(E)$, and the subdifferential of a function $f$ by $\partial f$.

\section{Setting of Problems}\label{sec:formulation}

\subsection{Maximum Hands-off Control}

We consider a continuous-time linear system defined by
\begin{equation}
\dot{x}(t)=Ax(t)+Bu(t), \quad 0 \leq t \leq T,
\label{eq:system}
\end{equation}
where 
$x(t)\in{\mathbb{R}}^n$ is the state vector,
$u(t)\in{\mathbb{R}}^m$ is the control input,
and $T>0$ is the final time of control. 
For the system~\eqref{eq:system}, 
a control input $u$ satisfying the following conditions is called {\em feasible}.
\begin{enumerate}
\item 
The state $x(t)$ is steered from a given initial state $x(0) = x_0\in\mathbb{R}^n$ 
to a given target state $x_f\in\mathbb{R}^n$ at time $T>0$ (i.e., $x(T)=x_f$ holds), and
\item
The magnitude constraint $\|u\|_{L^\infty} \leq 1$ holds.
\end{enumerate}
This paper assumes $x_f=0$ without losing generality.
The set of all feasible controls is denoted by $\mathcal{U}(x_0, T)$ (or simply $\mathcal{U}$), for given $x_0$ and $T$, i.e.,
\begin{equation*}
	\mathcal{U}(x_0, T) \triangleq \left\{u: \int_{0}^{T} e^{-At}Bu(t)dt = - x_0, \|u\|_{L^\infty}\leq 1 \right\}.
\end{equation*}
This paper considers a case where $\mathcal{U}(x_0, T)$ is not empty
(otherwise, it is obvious that the maximum hands-off control does not exist).
Here, the maximum hands-off control refers to a control $u\in\mathcal{U}$ having the minimum support length. 
In other words, this optimal control has the minimum $L^0$ cost among all control inputs in $\mathcal{U}$.
This problem is formulated as follows:

\begin{problem}[maximum hands-off control]\label{prob:main}
For given $x_0\in\mathbb{R}^n$ and $T>0$, find a control input $u$ on $[0, T]$ 
that minimizes $\|u\|_{L^0}$ subject to $u \in \mathcal{U}(x_0, T)$.
\end{problem}

Due to the $L^0$ norm in the cost function, it is difficult from the viewpoint of computational effort to solve Problem~\ref{prob:main} as it is.
Then, conventional studies adopt a convex approximation method to replace the $L^0$ norm by the $L^1$ norm, 
and evaluate the mathematical relationship between these two optimal controls.
In contrast, this paper considers a non-convex approximation method employing the non-convex penalty that can  better approximate the original $L^0$ norm.
Although the non-convex optimization problem is more difficult to analyze in comparison with the convex optimization problem, the non-convex optimization can potentially have the maximum hands-off control under a less restrictive condition. 
In fact, as proved in Theorem~\ref{thm:L0-L1}, this non-convex problem is always equivalent to the maximum hands-off control problem. 
Note that when specific two optimal control problems have the same solution set, those problems are said to be {\em equivalent}. 
Next, let us introduce the non-convex optimal control problems to discuss in this paper.

\subsection{Non-Convex Optimal Control Problems}
Based on the characteristics of non-convex penalty function that promotes sparsity 
such as the $\ell^p$ norm $(0<p<1)$, the MCP, and the SCAD, which are used in the signal/image processing fields,
this paper approximates the $L^0$ norm of $u$ by the difference in integral value between $\|u(t)\|_{\ell^1}$ and $\phi(u(t))$ over the interval $[0, T]$ with a given function $\phi:\mathbb{R}^m\to\mathbb{R}$.
This paper posits the following assumption on the function $\phi$.

\begin{assumption}\label{ass:phi}
The function $\phi$ satisfies the following:
\begin{enumerate}
\item[(A1)]\label{ass:phi_separable} $\phi$ is additively separable, i.e., 
	there exist functions $\phi_j: \mathbb{R}\to\mathbb{R}$, $j=1,2,\dots,m$, that satisfy
	$\phi(u)=\sum_{j=1}^{m} \phi_j(u_j)$,
	where $u = [u_1, u_2, \dots, u_m]^\top$.
\item[(A2)]\label{ass:phi_even} $\phi_j(u_j)=\phi_j(-u_j)$ on $[0, 1]$ for all $j=1,2,\dots,m$.
\item[(A3)]\label{ass:phi_bound} $\phi_j(0)=0$, 
			$\phi_j(u_j) < \phi_j(1) |u_j|$ on $(-1, 1) \backslash \{0\}$, 
			and $\phi_j(1) < 1$ for all $j=1,2,\dots,m$.
\item[(A4)]\label{ass:phi_terminal} $\phi_i(1) = \phi_j(1)$ for all $i,j=1,2,\dots,m$.
\end{enumerate}
\end{assumption}

With such function $\phi$, the non-convex optimal control problems to consider in this paper are formulated as follows:
\begin{problem}\label{prob:convex}
For given $x_0\in\mathbb{R}^n$ and $T>0$, find a control input $u$ on $[0, T]$ that solves the following:
\begin{equation*}
\begin{aligned}
  & \underset{u}{\text{minimize}}
  & & \|u\|_{L^1} - \int_{0}^{T} \phi(u(t)) dt\\
  & \text{subject to}
  & & u \in \mathcal{U}(x_0, T).
\end{aligned}
\end{equation*}
\end{problem}

Throughout the paper, the cost function of Problem~\ref{prob:convex} is denoted by $J:\mathcal{U}\to\mathbb{R}$, i.e., 
$J(u) \triangleq \|u\|_{L^1} - \int_{0}^{T} \phi(u(t)) dt$.

\begin{rmk}
This paper tries to derive optimal control problems that are equivalent to the maximum hands-off control problem.
Here, it is known that the maximum hands-off control is a bang-off-bang control~\cite{ITKKTAC18} (i.e., it takes only the discrete values belonging to the set $\{0, \pm 1\}^m$).
The assumptions~(A1), (A2), and (A3) are introduced to give the property of this to the optimal control (see the proof of Theorem~\ref{thm:discrete-L1} for details).
The assumption~(A4) is introduced to have the equivalence hold. 
If the assumption~(A4) does not hold, then the weights of sparsity varies according to the control variable $u_j$ in Problem~\ref{prob:convex}, and the equivalence does not necessarily hold (see the proof of Theorem~\ref{thm:L0-L1} for details).
In addition to the above-described technical reasons, when the assumption~(A3) holds, the integrand of the cost function takes its minimum value $0$ at its origin, i.e., $|u_j| - \phi_j(u_j)>0$ holds for any $u_j\in[-1,1]\backslash\{0\}$.
\end{rmk}

\begin{rmk}\label{rmk:penalties}
Examples of Problem~\ref{prob:convex} include the following optimal control problems.
Here, the cost function is denoted by 
$\int_{0}^{T} \sum_{j=1}^{m} \psi(u_j(t)) dt$
with a penalty function $\psi:\mathbb{R}\to\mathbb{R}$.
Then, the function $\phi$ is represented by 
$\phi(u(t)) = \sum_{j=1}^{m} \left( |u_j(t)| - \psi(u_j(t)) \right)$.
The functions $\psi$ and $\phi$ for each case are depicted in Fig.~\ref{fig:Lp-MCP-SCAD}.

\begin{itemize}
\item The $L^p$ penalty with $0<p<1$~\cite{FraFri93}:
\[
	\psi(u_j) = \lambda |u_j|^p
\]
where $\lambda>0$.

\item The Minimax Concave Penalty (MCP)~\cite{Zha10}:
\[
	\psi(u_j) = 
	\begin{cases}
		\lambda|u_j| - \frac{u_j^2}{2\alpha}, & \mbox{~if~} |u_j| \leq \alpha\lambda\\
		\frac{\alpha\lambda^2}{2}, & \mbox{~if~} |u_j| > \alpha\lambda
	\end{cases}
\]
where $\lambda>0$ and $\alpha>0$.

\item The Smoothly Clipped Absolute Deviation (SCAD)~\cite{FanLi01}:
\[
	\psi(u_j) = 
	\begin{cases}
		\lambda |u_j|, 
			& \mbox{~if~} |u_j| \leq \lambda\\
		-\frac{u_j^2-2\alpha\lambda|u_j|+\lambda^2}{2(\alpha-1)}, 
			& \mbox{~if~} \lambda < |u_j| \leq \alpha \lambda\\
		\frac{(\alpha+1)\lambda^2}{2}, 
			& \mbox{~if~} |u_j| > \alpha \lambda
	\end{cases}
\]
where $\lambda\in(0, 1)$ and $\alpha>1$.

\item The Log-Sum Penalty (LSP)~\cite{CanWakBoy08}:
\[
	\psi(u_j) = \lambda \log\left(1+\frac{|u_j|}{\alpha}\right)
\]
where $\lambda > 0$ and $\alpha>0$.

\item The capped $L^1$ penalty~\cite{Zhan10}:
\[
	\psi(u_j) = \lambda \min \left\{|u_j|, \alpha\right\}
\]
where $\lambda>0$ and $\alpha\in(0, 1)$.

\item The $L^1/L^2$ penalty:
\[
	\psi(u_j) = |u_j| - \lambda u_j^2
\]
where $\lambda\in(0,1)$ 
(i.e., $\psi$ is the difference between the $\ell^1$ norm and the squared $\ell^2$ norm of the control variables).
\end{itemize}

\begin{figure}[tb]
  \centering
    \includegraphics[width=0.47\linewidth]{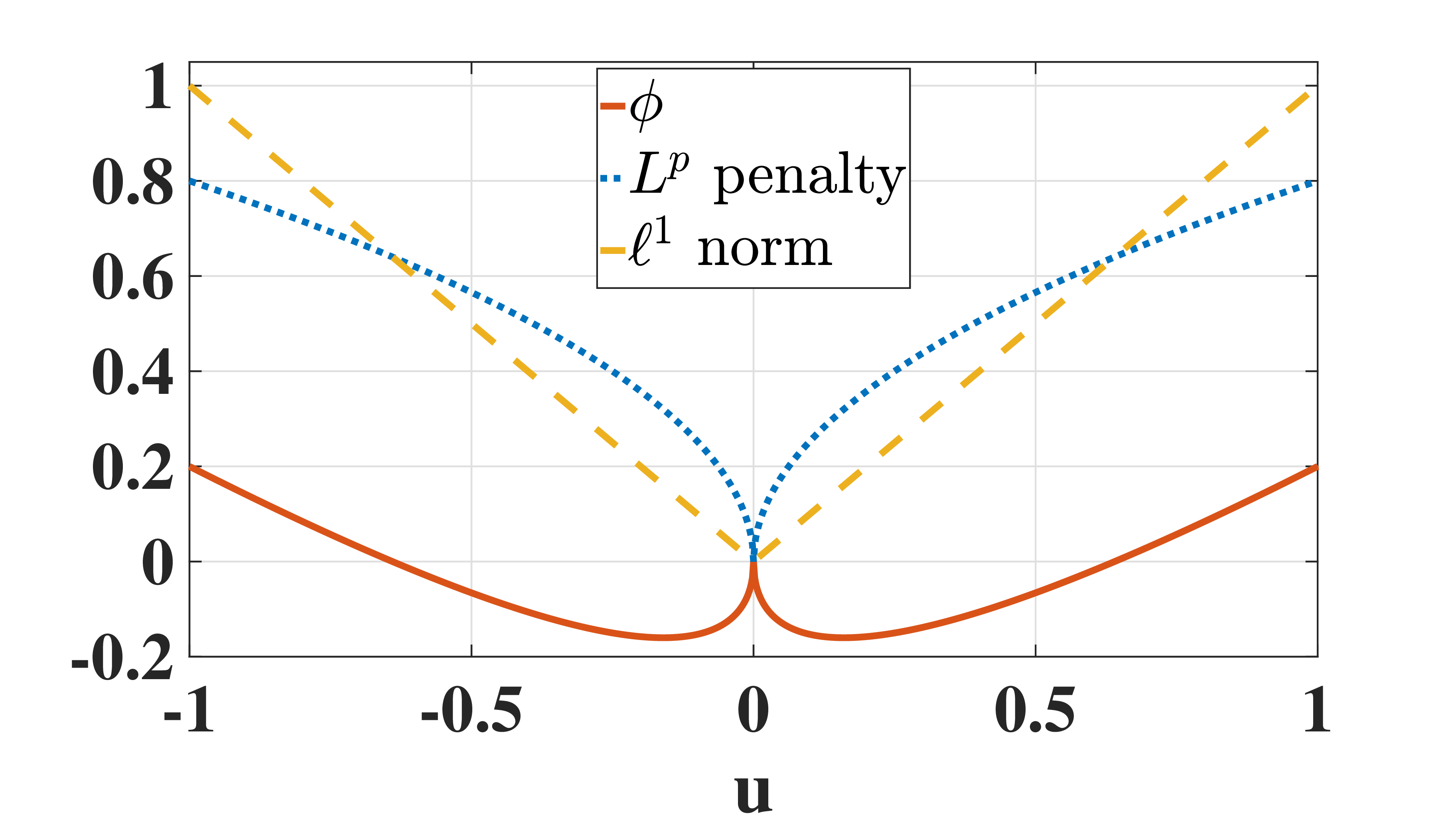}
	\includegraphics[width=0.47\linewidth]{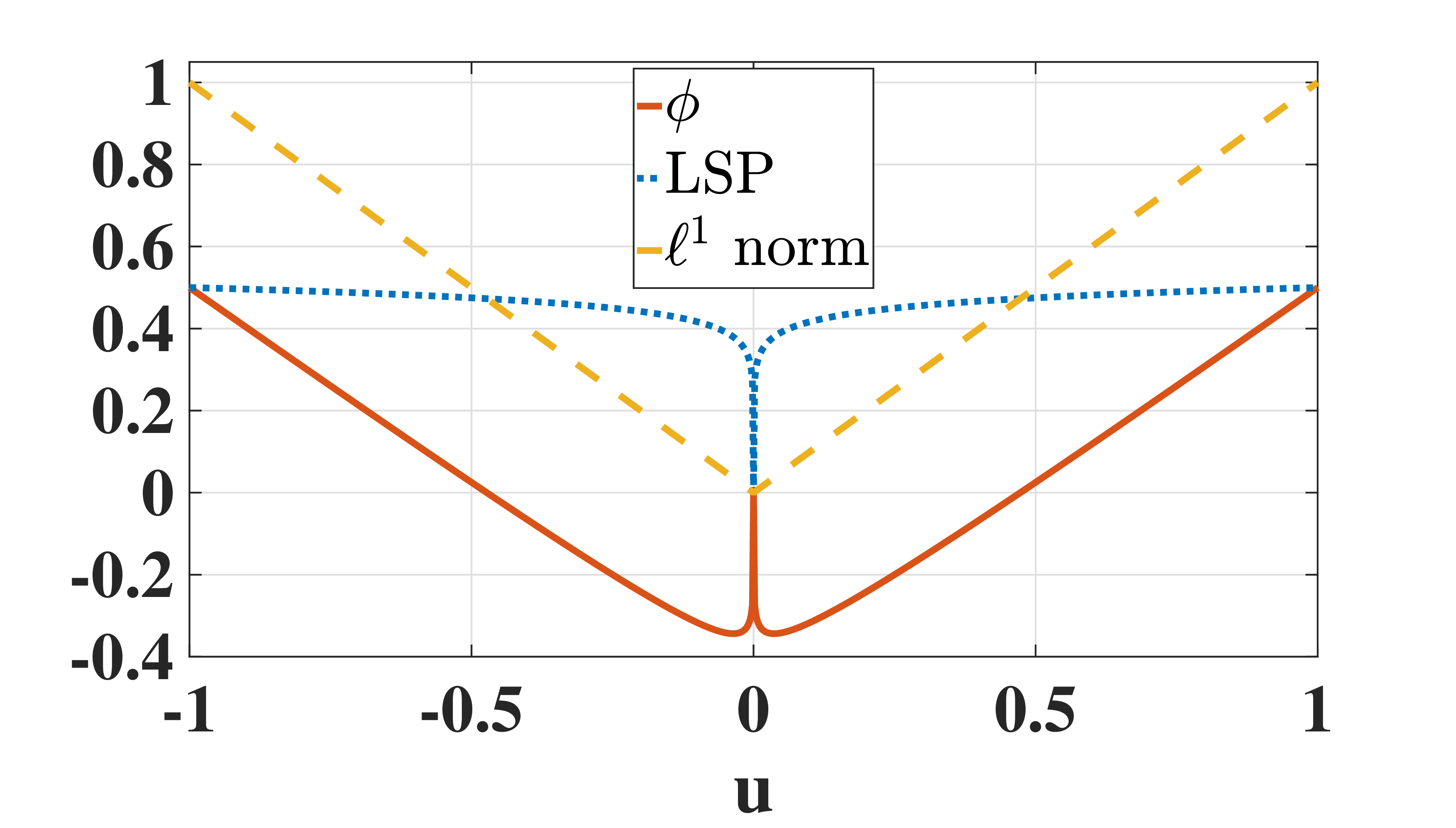}    
    \includegraphics[width=0.47\linewidth]{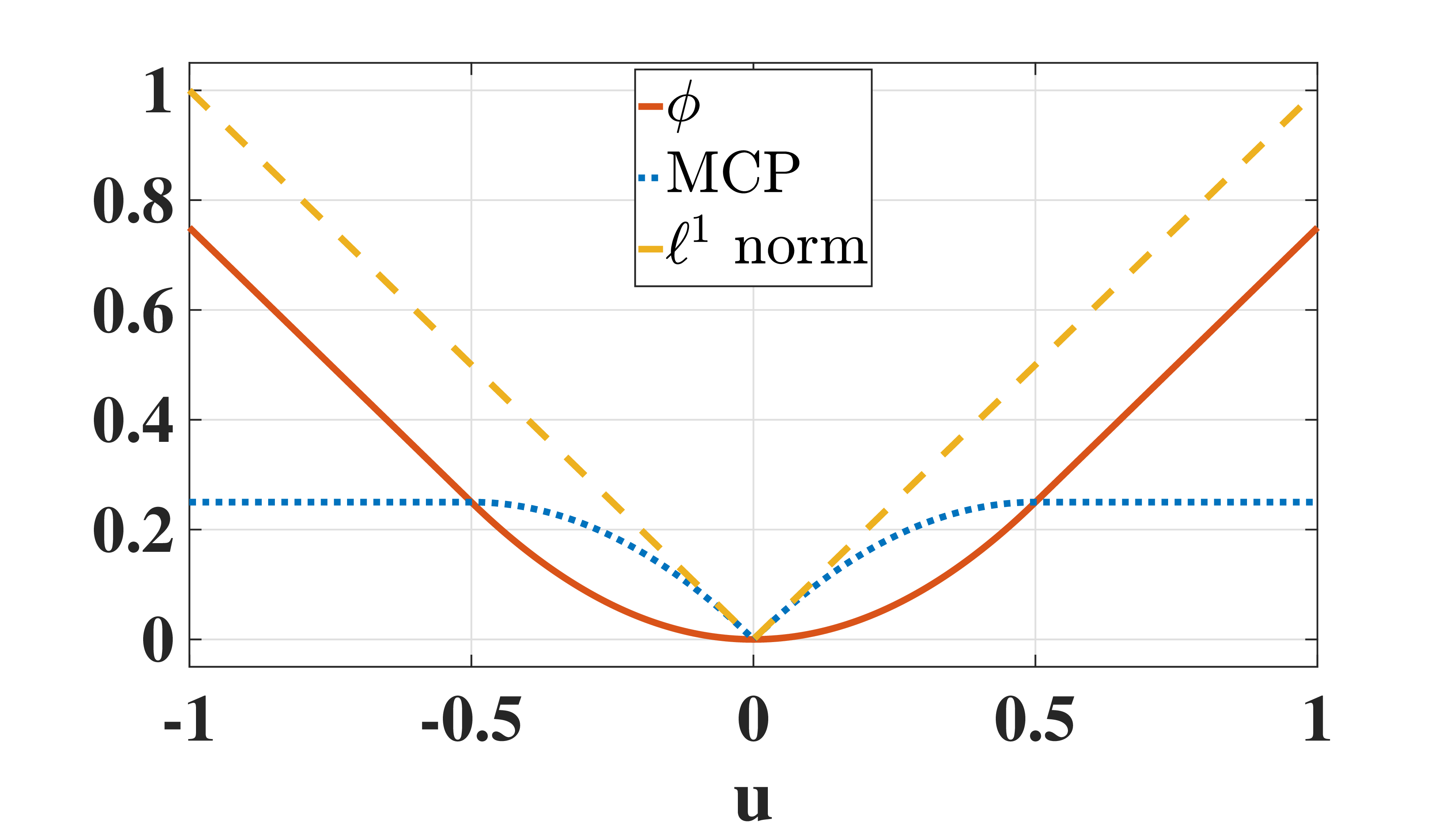}
    \includegraphics[width=0.47\linewidth]{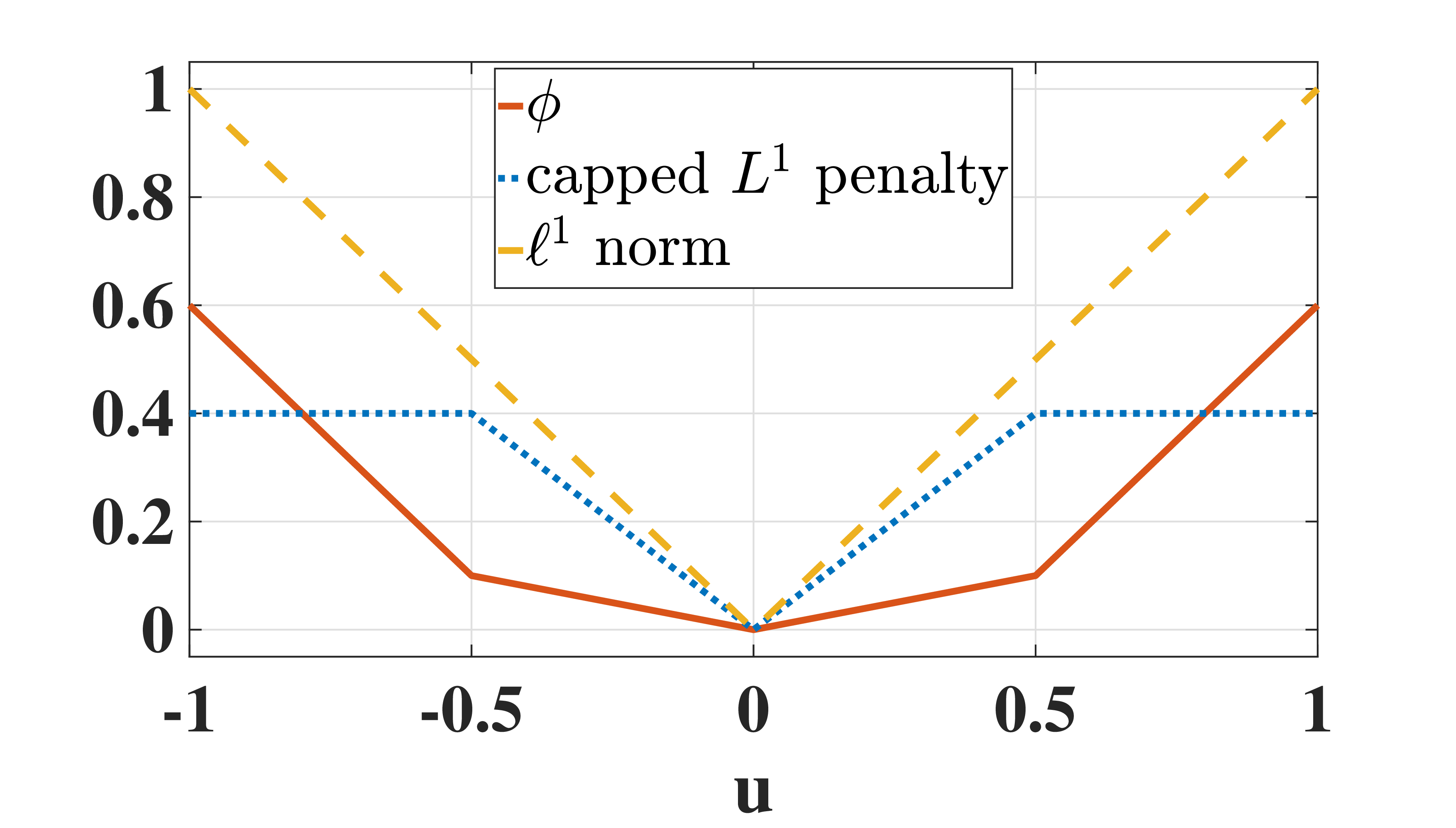}
    \includegraphics[width=0.47\linewidth]{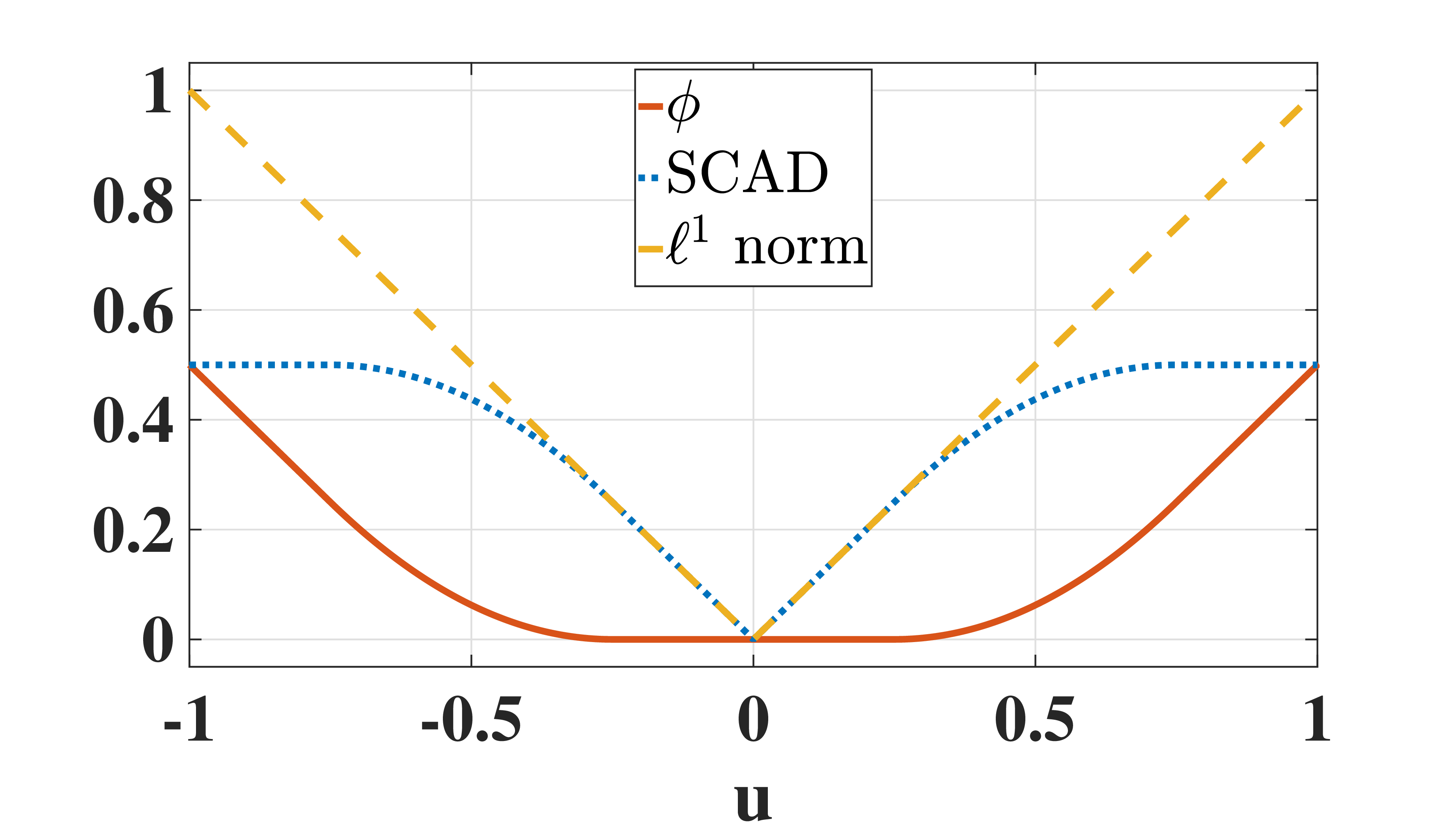}
    \includegraphics[width=0.47\linewidth]{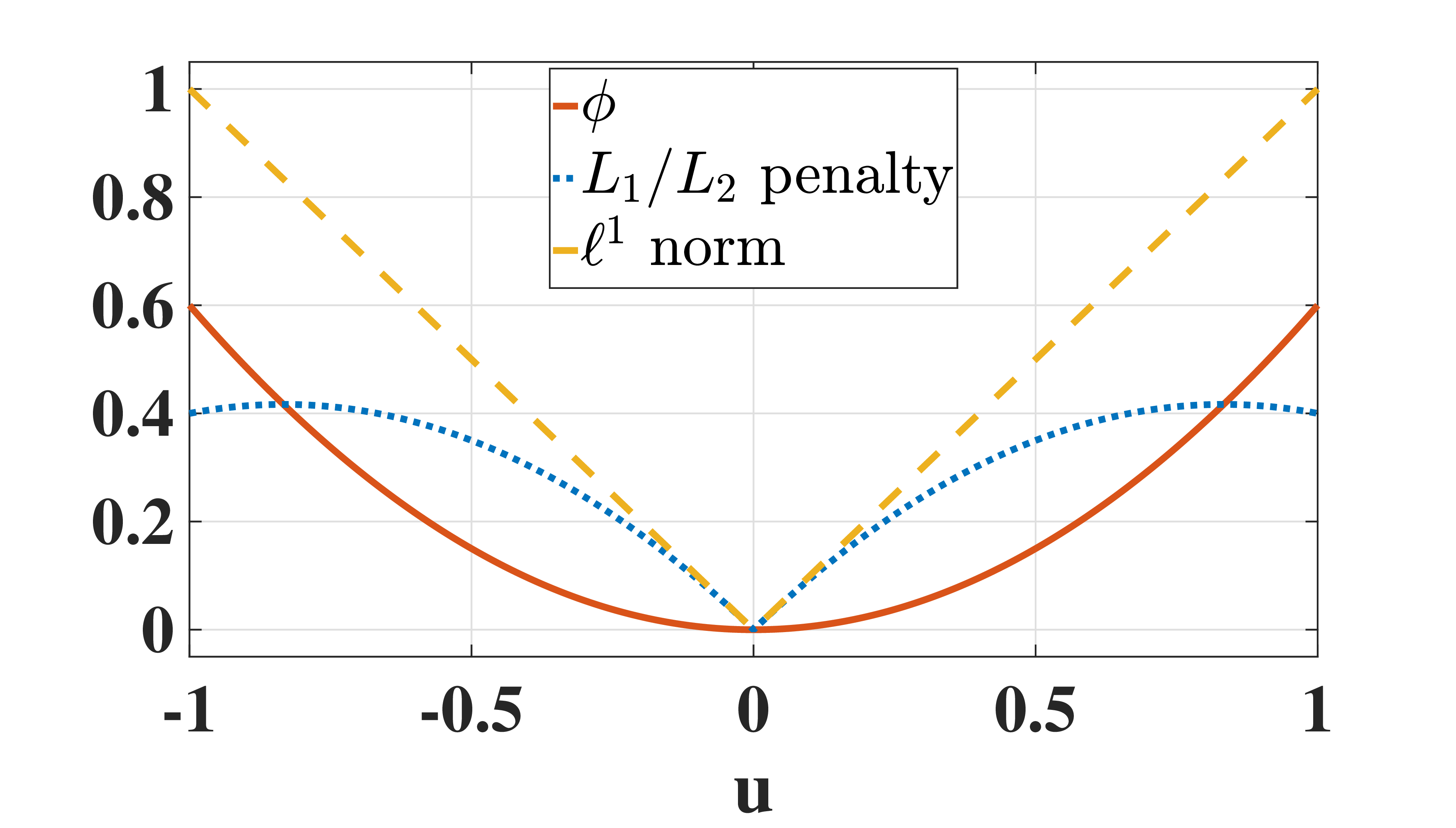}
  \caption{Examples of non-convex penalty $\psi$ (dotted line) and the function $\phi$ (solid line).
    The $L^p$ penalty with $(p,\lambda)=(0.5, 0.8)$ (top left);
    The MCP with $(\lambda, \alpha) = (0.25, 2)$ (middle left);
    The SCAD with $(\lambda, \alpha)=(0.25, 3)$ (bottom left);
    The LSP with $(\lambda, \alpha) = (0.5/\log(1+1/\alpha), 10^{-6})$ (top right);
    capped $L^1$ penalty with $(\lambda, \alpha) = (0.8, 0.5)$ (middle right);
    The $L^1/L^2$ penalty with $\lambda=0.6$ (bottom right).}
\label{fig:Lp-MCP-SCAD}
\end{figure}  
\end{rmk}

\section{Analysis and Computation}\label{sec:analysis_comp}

\subsection{Sparsity of the Non-convex Optimal Controls}
Here, the relationship between the maximum hands-off control (Problem~\ref{prob:main}) and the non-convex optimal controls (Problem~\ref{prob:convex}) is considered.
Firstly, the bang-off-bang property of the optimal solutions to Problem~\ref{prob:convex} is proved.
This property plays an important role in the theorem of equivalence (Theorem~\ref{thm:L0-L1}).
For this theorem, the following lemma is prepared.

\begin{lemma}[Theorem 8.2, \cite{HerLas}]\label{lem:existence_bang-off}
Let $\mathcal{I}\subset\mathbb{R}$ be any subset of the real line having finite Lebesgue measure,
\begin{align*}
	&\Psi\triangleq\{a \in L^{\infty}(\mathcal{I}): a(t) \in [0, 1]\},\\
	&\Psi^0\triangleq\{\chi_E: E \mbox{~a measurable subset of~} \mathcal{I}\},
\end{align*}	
and $y$ be a function with components $y_1, y_2, \dots, y_n \in L^1({\mathcal{I}})$,
where $\chi_E$ is the characteristic function of $E$, 
i.e., $\chi_E(t)=1$ for $t\in E$ and $\chi_E(t)=0$ for $t\not\in E$.
Then, 
\[
	\left\{ \int_{\mathcal{I}} y(t) a(t) dt : a \in \Psi \right\}
	= \left\{ \int_{\mathcal{I}} y(t) a(t) dt : a \in \Psi^0 \right\}.
\]
\end{lemma}

The following theorem guarantees the bang-off-bang property of the non-convex optimal control.
Although this theorem assumes the existence of the optimal solution, it should be noted that the existence will be proved later in Theorem~\ref{thm:L0-L1}.

\begin{theorem}[bang-off-bang property]\label{thm:discrete-L1}
Suppose that the function $\phi$ satisfies (A1), (A2), and (A3) in Assumption~\ref{ass:phi}.
Then, any optimal solution to Problem~\ref{prob:convex} takes values belonging to the set $\{0, \pm 1\}^m$ almost everywhere.
\end{theorem}
\begin{IEEEproof}
Here, it is assumed that the optimal solution to Problem~\ref{prob:convex} exists,
and let us take any optimal solution $\hat{u}\in\mathcal{U}$.
To show this bang-off-bang property, let us suppose that
\begin{equation}\label{eq:ass_discrete}
	\mu \left( \left\{ t\in[0, T]: \hat{u}_j(t) \not\in \{0, \pm1\} \right\} \right) > 0
\end{equation}
holds for some $j$ and show this leads to a contradiction.
Put 
\begin{align*}
	&\mathcal{I}_{+}^{(1)} \triangleq \left\{t \in [0, T]: \hat{u}_j(t) = 1 \right\}, \\
	&\mathcal{I}_{+}^{(2)} \triangleq \left\{t \in [0, T]: \hat{u}_j(t) \in (0, 1) \right\}, \\
	&\mathcal{I}_{0} \triangleq \left\{t \in [0, T]: \hat{u}_j(t) = 0 \right\}, \\
	&\mathcal{I}_{-}^{(1)} \triangleq \left\{t \in [0, T]: \hat{u}_j(t) = -1 \right\},\\
	&\mathcal{I}_{-}^{(2)} \triangleq \left\{t \in [0, T]: \hat{u}_j(t) \in (-1, 0) \right\}.
\end{align*}
These sets are mutually disjoint, and their union is the interval $[0, T]$.
From~\eqref{eq:ass_discrete}, note that
\begin{align}\label{eq:ass_discrete_2}
	\mu(\mathcal{I}_{+}^{(2)}) + \mu(\mathcal{I}_{-}^{(2)}) > 0.
\end{align}
Note also that 
\begin{align}\label{eq:feasibility}
	-x_{0} = \xi_{+}^{(1)} + \xi_{+}^{(2)} + \xi_{-}^{(1)} + \xi_{-}^{(2)} 
	+ \sum_{i\neq j}  \int_{0}^{T} e^{-At} b_i \hat{u}_i(t) dt,
\end{align}
where $b_i\in\mathbb{R}^n$ is the $i$th column of the matrix $B$, and
\begin{align*}
	\xi_{+}^{(1)} \triangleq \int_{\mathcal{I}_{+}^{(1)}} e^{-At} b_j \hat{u}_j(t) dt,\quad
	\xi_{+}^{(2)} \triangleq \int_{\mathcal{I}_{+}^{(2)}} e^{-At} b_j \hat{u}_j(t) dt,\\
	\xi_{-}^{(1)} \triangleq \int_{\mathcal{I}_{-}^{(1)}} e^{-At} b_j \hat{u}_j(t) dt,\quad
	\xi_{-}^{(2)} \triangleq \int_{\mathcal{I}_{-}^{(2)}} e^{-At} b_j \hat{u}_j(t) dt.
\end{align*}
Here, from Lemma~\ref{lem:existence_bang-off},
there exist functions $w_{+}^{(2)}$ and $w_{-}^{(2)}$ satisfying 
\begin{align}
	&w_{+}^{(2)}(t)\in
		\begin{cases} 
			\{0, 1\}, &\forall t\in\mathcal{I}_{+}^{(2)},\\
			\{0\}, &\forall t\in [0,T] \backslash \mathcal{I}_{+}^{(2)},
		\end{cases} \label{eq:w_p_01}\\
	&\int_{\mathcal{I}_{+}^{(2)}} 
		\begin{bmatrix}e^{-At} b_j \\ 1 \end{bmatrix} \hat{u}_{j}(t) dt
		 = \int_{\mathcal{I}_{+}^{(2)}} 
		\begin{bmatrix}e^{-At} b_j \\ 1 \end{bmatrix} w_{+}^{(2)}(t) dt, \label{eq:w_p_feasibility}\\
	&w_{-}^{(2)}(t)\in
		\begin{cases} 
			\{0, -1\}, &\forall t\in\mathcal{I}_{-}^{(2)},\\
			\{0\}, &\forall t\in [0,T] \backslash \mathcal{I}_{-}^{(2)},
		\end{cases}\label{eq:w_m_01}\\
	&\int_{\mathcal{I}_{-}^{(2)}} 
		\begin{bmatrix}e^{-At} b_j \\ 1 \end{bmatrix} \hat{u}_{j}(t) dt
		 = \int_{\mathcal{I}_{-}^{(2)}} 
		\begin{bmatrix}e^{-At} b_j \\ 1 \end{bmatrix} w_{-}^{(2)}(t) dt. \label{eq:w_m_feasibility}
\end{align}
Define 
\[
	w_j(t) \triangleq
	\begin{cases} 
		\hat{u}_{j}(t), &\forall t\in\mathcal{I}_{+}^{(1)}\cup\mathcal{I}_{0}\cup\mathcal{I}_{-}^{(1)},\\
		w_{+}^{(2)}(t), &\forall t\in\mathcal{I}_{+}^{(2)},\\
		w_{-}^{(2)}(t), &\forall t\in\mathcal{I}_{-}^{(2)},
	\end{cases}
\]
and $w(t)\triangleq[\hat{u}_1(t), \dots, \hat{u}_{j-1}(t), w_j(t), \hat{u}_{j+1}(t),\dots, \hat{u}_m(t)]^\top$ on $[0, T]$.
Note that $w_j(t)\in\{0, \pm 1\}$ on $[0,T]$ from~\eqref{eq:w_p_01}, \eqref{eq:w_m_01}, and the definition of the sets $\mathcal{I}_{+}^{(1)}$, $\mathcal{I}_{0}$, and $\mathcal{I}_{-}^{(1)}$.
Note also that $w\in\mathcal{U}$ from~\eqref{eq:feasibility}, \eqref{eq:w_p_feasibility}, and \eqref{eq:w_m_feasibility}.
Furthermore,
\begin{align}\label{eq:w_L1_norm}
\begin{split}
	&\int_{0}^{T} |w_j(t)| dt \\
	&= \int_{\mathcal{I}_{+}^{(1)} \cup \mathcal{I}_{0} \cup \mathcal{I}_{-}^{(1)}}|\hat{u}_j(t)| dt 
	+ \int_{\mathcal{I}_{+}^{(2)}} w_j(t) dt - \int_{\mathcal{I}_{-}^{(2)}} w_j(t) dt \\
	&= \int_{\mathcal{I}_{+}^{(1)} \cup \mathcal{I}_{0} \cup \mathcal{I}_{-}^{(1)}}|\hat{u}_j(t)| dt 
	+ \int_{\mathcal{I}_{+}^{(2)}} \hat{u}_j(t) dt - \int_{\mathcal{I}_{-}^{(2)}} \hat{u}_j(t) dt \\
	&= \int_{0}^{T} |\hat{u}_j(t)| dt,
\end{split}
\end{align}
where the second relation follows from~\eqref{eq:w_p_feasibility} and~\eqref{eq:w_m_feasibility}.
Here, from~\eqref{eq:ass_discrete_2}, $\mu(\mathcal{I}_{+}^{(2)}) > 0$ or $\mu(\mathcal{I}_{-}^{(2)}) > 0$ holds.
When $\mu(\mathcal{I}_{+}^{(2)}) > 0$,
\begin{align*}
	\int_{\mathcal{I}_{+}^{(2)}} \phi_j(\hat{u}_j(t)) dt 
	& <  \phi_j(1) \int_{\mathcal{I}_{+}^{(2)}} \hat{u}_j(t) dt \\
	& =  \phi_j(1) \int_{\mathcal{I}_{+}^{(2)}} w_j(t) dt \\
	& =  \int_{\mathcal{I}_{+}^{(2)}} \phi_j(w_j(t)) dt
\end{align*}
holds.
The first relation follows from the assumption~(A3), 
the second relation from~\eqref{eq:w_p_feasibility},
and the third relation from $w_j(t)\in\{0, 1\}$ on $\mathcal{I}_{+}^{(2)}$ and the assumption~(A3).
When $\mu(\mathcal{I}_{+}^{(2)}) = 0$, 
$\int_{\mathcal{I}_{+}^{(2)}} \phi_j(\hat{u}_j(t)) dt 
= \int_{\mathcal{I}_{+}^{(2)}} \phi_j(w_j(t)) dt
= 0$
holds.
In the same way, when $\mu(\mathcal{I}_{-}^{(2)}) > 0$,
\begin{align*}
	\int_{\mathcal{I}_{-}^{(2)}} \phi_j(\hat{u}_j(t)) dt <  \int_{\mathcal{I}_{-}^{(2)}} \phi_j(w_j(t)) dt
\end{align*}
holds from the assumption~(A2).
When $\mu(\mathcal{I}_{-}^{(2)}) = 0$, 
$\int_{\mathcal{I}_{-}^{(2)}} \phi_j(\hat{u}_j(t)) dt 
= \int_{\mathcal{I}_{-}^{(2)}} \phi_j(w_j(t)) dt
= 0$
holds.
From these,
\begin{align}\label{eq:w_phi}
\begin{split}
	&\int_{0}^{T} \phi_j(\hat{u}_j(t)) dt \\
	&= \int_{\mathcal{I}_{+}^{(1)} \cup \mathcal{I}_{0} \cup \mathcal{I}_{-}^{(1)}} \phi_j(\hat{u}_j(t)) dt 
	+ \int_{\mathcal{I}_{+}^{(2)} \cup \mathcal{I}_{-}^{(2)}} \phi_j(\hat{u}_j(t)) dt \\
	&< \int_{\mathcal{I}_{+}^{(1)} \cup \mathcal{I}_{0} \cup \mathcal{I}_{-}^{(1)}} \phi_j(\hat{u}_j(t)) dt 
	+ \int_{\mathcal{I}_{+}^{(2)} \cup \mathcal{I}_{-}^{(2)}} \phi_j(w_j(t)) dt \\
	&= \int_{0}^{T} \phi_j(w_j(t)) dt.
\end{split}
\end{align}
Therefore, from~\eqref{eq:w_L1_norm}, \eqref{eq:w_phi}, and the assumption~(A1), 
\begin{align*}
	J(\hat{u})
	&= \int_{0}^{T} \sum_{i\neq j} \left(|\hat{u}_i(t)| - \phi_i(\hat{u}_i(t)) \right) dt\\
	&\mspace{30mu} + \int_{0}^{T} \left(|\hat{u}_j(t)| - \phi_j(\hat{u}_j(t)) \right) dt\\
	&> \int_{0}^{T} \sum_{i\neq j} \left(|\hat{u}_i(t)| - \phi_i(\hat{u}_i(t)) \right) dt\\
	&\mspace{30mu} + \int_{0}^{T} \left(|w_j(t)| - \phi_j(w_j(t)) \right) dt\\
	&= J(w)
\end{align*}
holds in contradiction to the optimality of $\hat{u}$.
This completes the proof.
\end{IEEEproof}

Next, the equivalence between Problem~\ref{prob:main} and Problem~\ref{prob:convex}, the main result of this paper, is proved.
For this purpose, we recall a property of the maximum hands-off control.

\begin{lemma}[Theorem 3, \cite{ITKKTAC18}]\label{lem:recall_L0}
Problem~\ref{prob:main} is considered.
It is supposed that the set $\mathcal{U}(x_0, T)$ is not empty.
Then, there exists at least one optimal solution, 
and any optimal solution takes values belonging to the set $\{0, \pm 1\}^m$ almost everywhere.
\end{lemma}

\begin{theorem}[existence and equivalence]
\label{thm:L0-L1}
It is supposed that the function $\phi$ satisfies Assumption~\ref{ass:phi} and the set $\mathcal{U}(x_0, T)$ is not empty.
Let ${\mathcal{U}}_{1}^{\ast}$ and ${\mathcal{U}}_{2}^{\ast}$ denote the set of all optimal solutions to Problem~\ref{prob:main} and Problem~\ref{prob:convex}, respectively.
Then, the sets ${\mathcal{U}}_{1}^{\ast}$ and ${\mathcal{U}}_{2}^{\ast}$ are not empty, 
and ${\mathcal{U}}_{1}^{\ast}={\mathcal{U}}_{2}^{\ast}$ holds.
\end{theorem}
\begin{IEEEproof}
The set $\mathcal{U}_1^\ast$ is not empty from Lemma~\ref{lem:recall_L0}.
From Theorem~\ref{thm:discrete-L1}, any optimal solution to Problem~\ref{prob:convex} takes  values belonging to $\{0, \pm 1\}^m$ almost everywhere.
Therefore, Problem~\ref{prob:convex} is rewritten as the optimal control problem to minimize $J(u)$ subject to $u \in \mathcal{U}'$, where 
\begin{align*}
	\mathcal{U}' \triangleq 
	\Bigg\{u: 
	&\int_{0}^{T} e^{-At}Bu(t)dt = -x_0, \\
	&u(t)\in\{0, \pm 1\}^m \mbox{~for~almost~all~} t \in [0, T]\Bigg\}.
\end{align*}
For any $u\in\mathcal{U}'$,
\begin{align}\label{eq:J_BOB}
\begin{split}
	J(u) 
	&= \sum_{j=1}^{m} \int_{0}^{T} \left( |u_j(t)| - \phi_j(u_j(t)) \right) dt \\	
	&= \sum_{j=1}^{m} \left( 1 - \phi_j(1) \right)\|u_j\|_{L^0}\\
	&= c \|u\|_{L^0}
\end{split}
\end{align}
holds, where the third equation follows from the assumption~(A4), 
and $c\triangleq 1 - \phi_j(1)>0$ for all $j$ from the assumption~(A3).
Therefore, from Lemma~\ref{lem:recall_L0}, any $\tilde{u}\in\mathcal{U}_1^\ast$ satisfies
\begin{equation*}\label{eq:U2_sub_U1}
	J(\tilde{u}) 
	= c \|\tilde{u}\|_{L^0}
	\leq c \|u\|_{L^0}
	= J(u)
\end{equation*}
for all $u\in\mathcal{U}'$.
This means $\tilde{u}\in\mathcal{U}_2^\ast$.
Therefore, $\mathcal{U}_1^\ast \subset \mathcal{U}_2^\ast$ holds,
and the set $\mathcal{U}_2^\ast$ is not empty.

Next, let us take any $\hat{u}\in\mathcal{U}_2^\ast$.
From $\hat{u}\in\mathcal{U}'$, 
\[
	\|\hat{u}\|_{L^0} 
	= \frac{1}{c} J(\hat{u})
	= \frac{1}{c} J(\tilde{u})
	= \|\tilde{u}\|_{L^0} 
\]
holds, where the first and the third equations follow from~\eqref{eq:J_BOB},
and the second equation follows from $\tilde{u} \in \mathcal{U}_2^\ast$.
This means $\hat{u}\in\mathcal{U}_1^\ast$, and $\mathcal{U}_2^\ast \subset \mathcal{U}_1^\ast$ holds.
\end{IEEEproof}

\begin{rmk}
In the previous study~\cite{NagQueNes16},
the conditions of making the maximum hands-off control problem and the $L^1$ optimal control problem equivalent to each other are analyzed using Pontryagin's maximum principle.
According to the results produced by the convex approximation method, when the matrix $A$ is non-singular and the system $(A, b_j)$ is controllable for all $j$, these two problems become equivalent.
In comparison, in the non-convex optimal control problem proposed by this paper, even if these conditions regarding the system are not satisfied, the equivalence to the maximum hands-off control problem always holds. 
This property will be confirmed with a numerical example in Section~\ref{sec:example}. 
\end{rmk}

\subsection{Numerical Optimization}\label{sec:computation}

Here, a numerical algorithm for the maximum hands-off control is provided based on Theorem~\ref{thm:L0-L1}.
To this end, let us reformulate Problem~\ref{prob:convex}.

\begin{proposition}\label{prop:convex_reformulation}
It is supposed that the function $\phi$ satisfies Assumption~\ref{ass:phi} and the set $\mathcal{U}(x_0, T)$ is not empty.
Let us define the following optimal control problem:
\begin{equation}\label{prob:convex_2}
\begin{aligned}
  & \underset{v, w}{\text{minimize}}
  & & J(v) + J(w) \\
  & \text{subject to}
  & & \dot{x}(t) = Ax(t) + \begin{bmatrix}B & -B\end{bmatrix} \begin{bmatrix} v(t) \\ w(t) \end{bmatrix},\\
  & & & x(0) = x_0, \quad x(T) = 0, \\
  & & & v(t), w(t) \in [0, 1]^m \mbox{~for~almost~all~} t\in[0,T].
\end{aligned}
\end{equation}
Then, the following holds.
\begin{enumerate}
\item[(i)] For any optimal solution $u^\ast$ to Problem~\ref{prob:convex}, define 
	\begin{equation}\label{eq:u_v_w_P2}
		v^\ast(t) \triangleq \max \{u^\ast(t), 0\},\quad
		w^\ast(t) \triangleq \max \{-u^\ast(t), 0\}
	\end{equation}
	where $t\in[0, T]$.
	Then, $(v^\ast, w^\ast)$ is the optimal solution to the problem~\eqref{prob:convex_2}.
\item[(ii)] For any optimal solution $(\hat{v}, \hat{w})$ to the problem~\eqref{prob:convex_2}, define 
	$\hat{u}(t) \triangleq \hat{v}(t) - \hat{w}(t)$, where $t\in[0, T]$.
	Then, $\hat{u}$ is the optimal solution to Problem~\ref{prob:convex}.
\end{enumerate}
\end{proposition}
\begin{IEEEproof}
See Appendix.
\end{IEEEproof}

The numerical method for the problem~\eqref{prob:convex_2} is provided using a time discretization approach~\cite[Section 2.3]{Ste}.
Firstly, the interval $[0,T]$ is divided into $N$ subintervals: $[0,T] = [0,\Delta) \cup \dots \cup [(N-1)\Delta, N\Delta]$,
where $\Delta$ is the discretized step size satisfying $T=N \Delta$.
Approximation is made with the state $x$ and the control $(v, w)$ as constants over each subinterval.
Then, for each $t=0, \Delta, \dots, N\Delta$,
the continuous-time system in the problem~\eqref{prob:convex_2} is described as
\[
 x_d[k+1] = \mathcal{A} x_d[k] + \mathcal{B} \begin{bmatrix} v_d[k] \\ w_d[k] \end{bmatrix}, \quad k = 0, 1, \dots, N-1
\]
where $x_d[k]\triangleq x(k\Delta)$, $v_d[k]\triangleq v(k\Delta)$,  $w_d[k]\triangleq w(k\Delta)$, and
\[
  \mathcal{A} \triangleq e^{A\Delta},\quad
  \mathcal{B} \triangleq\ \int_0^\Delta e^{At} \begin{bmatrix} B & -B \end{bmatrix} dt.
\]
The vector composed of the control variables is defined as
\[
  z \triangleq 
  \begin{bmatrix} v_d[0]^{\top}, w_d[0]^\top, \dots, v_d[N-1]^\top, w_d[N-1]^\top \end{bmatrix} ^\top.
\] 
Then, the state constraint $x(T)=0$ is approximated by
\[
 0 = x_d[N] = \zeta + \Phi z,
\]
where $\zeta \triangleq \mathcal{A}^N x_0$, and
\[
	\Phi \triangleq \begin{bmatrix} 
 		\mathcal{A}^{N-1}\mathcal{B} & \mathcal{A}^{N-2}\mathcal{B} & \dots & \mathcal{B}
 	\end{bmatrix} \in{\mathbb{R}}^{n\times mN}.
\]
Thus, the optimal control problem is approximated as
\begin{equation}
 \begin{aligned}
  & \underset{z}{\text{minimize}}
  & & J_d(z)\\
  & \text{subject to}
  & & z \in \mathcal{U}_d
 \end{aligned}
 \label{prob:DC}
\end{equation}
where
\begin{align*}
	& J_d(z) \triangleq g(z) - h(z),\\
	& g(z)\triangleq \|z\|_{\ell^1},\quad 
		h(z) \triangleq \sum_{k=0}^{N-1} \left(\phi(v_d[k]) + \phi(w_d[k]) \right),\\
	& \mathcal{U}_d \triangleq \{z\in\mathbb{R}^{2mN}: z \in [0, 1]^{2mN}, ~ \Phi z + \zeta = 0\}.
\end{align*}

\begin{algorithm} [t]
\caption{Numerical computation method for the maximum hands-off control based on the DC representation~\eqref{prob:DC}}\label{alg:DCA}
\begin{algorithmic}[1]
	\STATE Choose $z[0]\in\mathbb{R}^{2mN}$.
	\FOR {$l=0,1,2,\dots,$}
		\STATE $s[l] \in \partial h(z[l])$
		\STATE $z[l+1] \in \argmin_{z\in\mathcal{U}_d} (g(z)-s[l]^\top z)$
	\ENDFOR  
	\STATE Denote the found numerical solution by 
	\[
		z^\ast = 
		\begin{bmatrix} 
		v_d^\ast[0]^{\top}, w_d^\ast[0]^\top, \dots, v_d^\ast[N-1]^\top, w_d^\ast[N-1]^\top 
		\end{bmatrix}^\top,
		\]
		and compute $u_d^\ast[k]=v_d^\ast[k]-w_d^\ast[k]$ for all $k$.
	\STATE Employ the control 
	$u^\ast(t)=u_d^\ast[k]$ on each subinterval $[k\Delta, (k+1)\Delta)$ 
	as the numerical solution to Problem~\ref{prob:main}.
	\end{algorithmic}
\end{algorithm}

Thus, Problem~\ref{prob:convex} reduces to the problem~\eqref{prob:DC}.
Since the cost function $g - h$ is not convex, it is generally difficult to find a global optimal solution to this problem.
On the other hand, in the sparse optimization field, the function $\phi$ is convex on some area.
Taking this into consideration, $\phi$ is confined to be convex over $[0,1]^m$ in what follows.
(Each penalty introduced in Remark~\ref{rmk:penalties} certainly satisfies this property.
Note here that the convexity is not assumed over $[-1,1]^m$ to include the $L^p$ penalty and the LSP in the target.)
Then, the functions $g$ and $h$ are convex, and the problem~\eqref{prob:DC} belongs to the class called the {\em difference of convex functions (DC) optimization problem}.
Such a problem has been actively studied, and numerical algorithms to find a stationary point have been proposed~\cite{LePha18dc,TaoAn97,SunYinCheJia18}. 
This paper uses a method quoting the best-known DC algorithm (see Algorithm~\ref{alg:DCA}).
For the convex optimization problem in step~4,
the alternating direction method of multipliers~\cite{ADMMBoyd}
or numerical software packages such as CVX of MATLAB~\cite{cvx} can be used.

\section{Example}
\label{sec:example}

\begin{figure}[tb]
  \centering
    \includegraphics[width=0.47\linewidth]{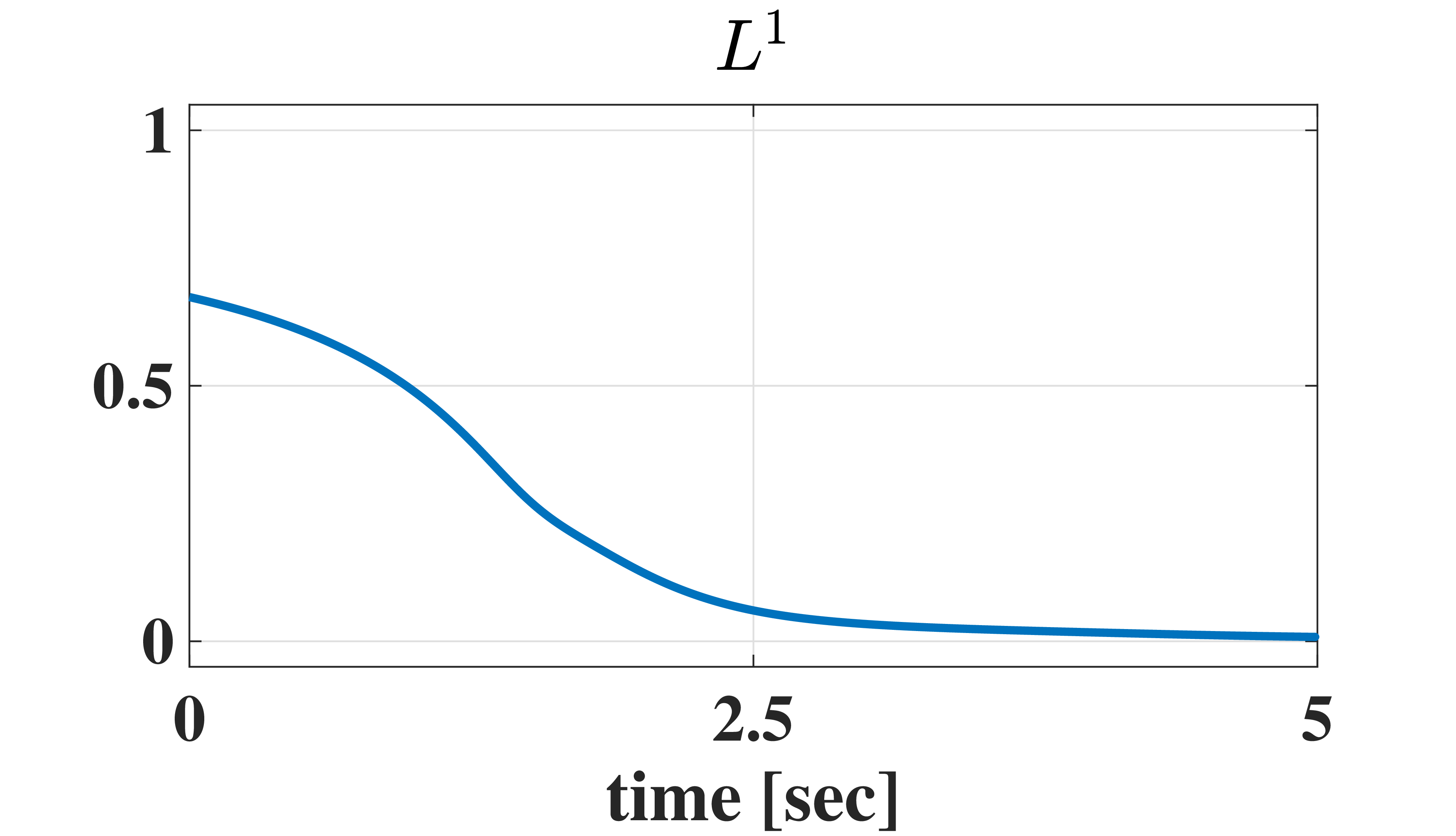}
    \includegraphics[width=0.47\linewidth]{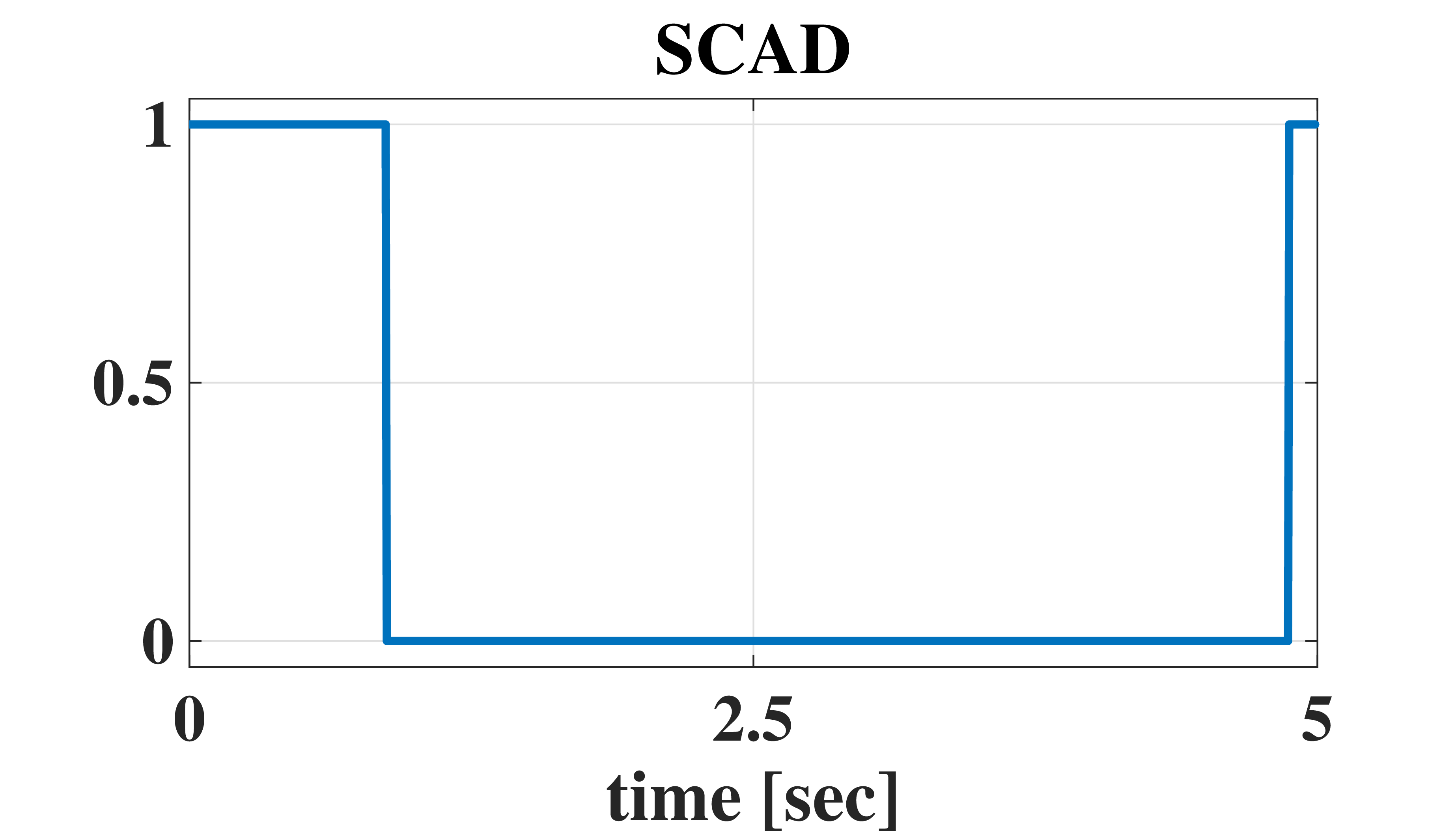}
    \includegraphics[width=0.47\linewidth]{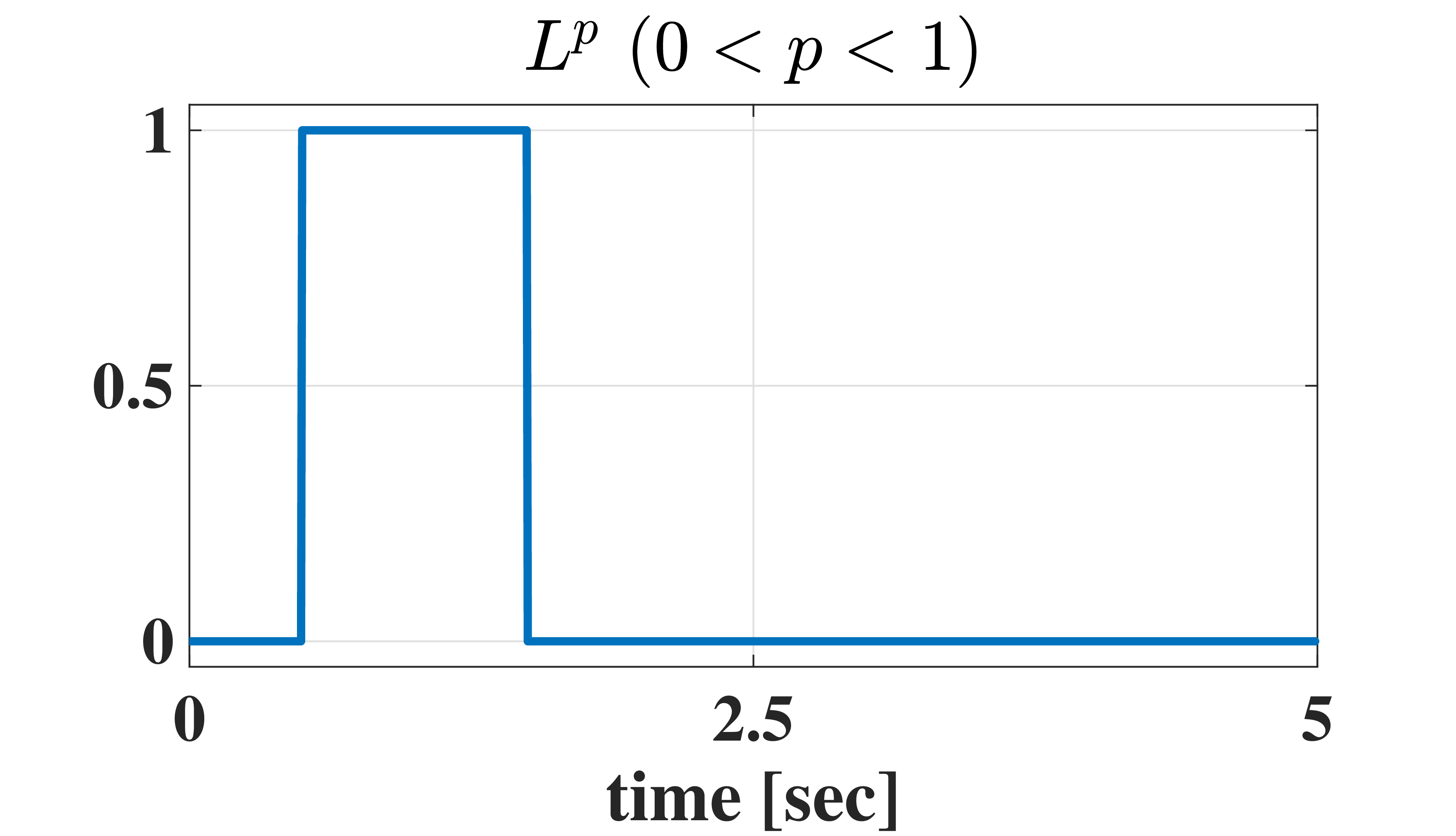}
    \includegraphics[width=0.47\linewidth]{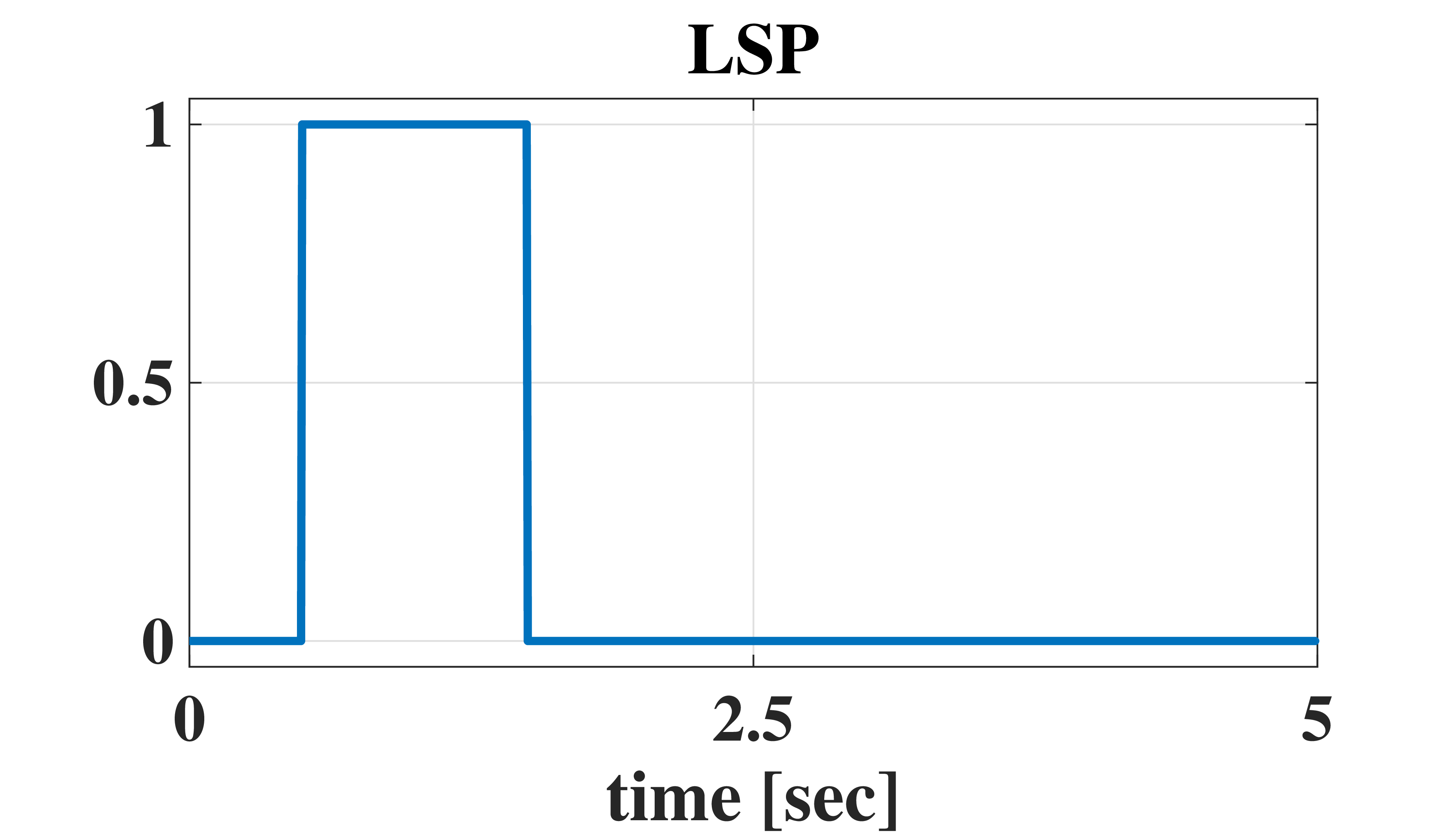}
    \includegraphics[width=0.47\linewidth]{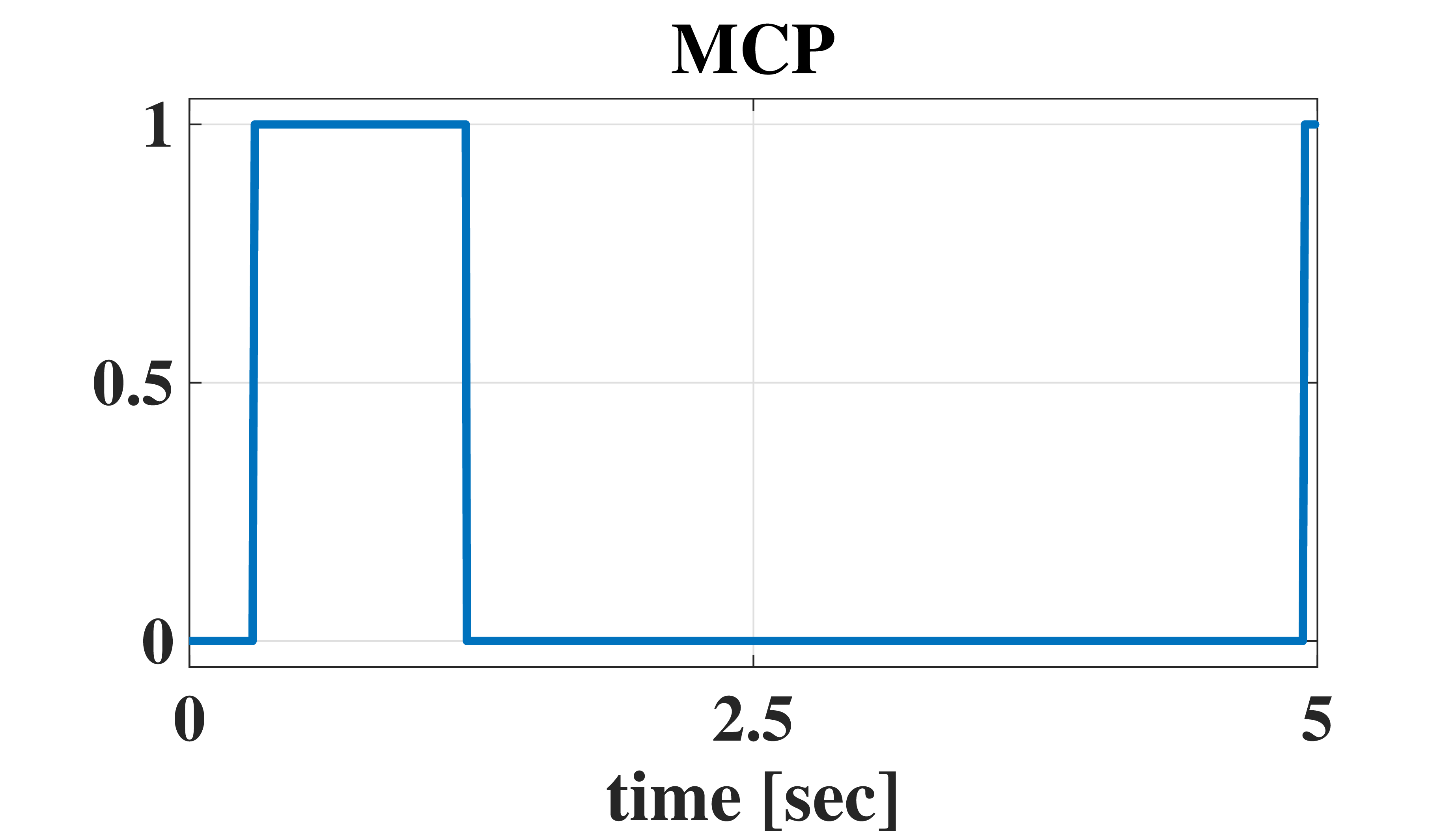}
    \includegraphics[width=0.47\linewidth]{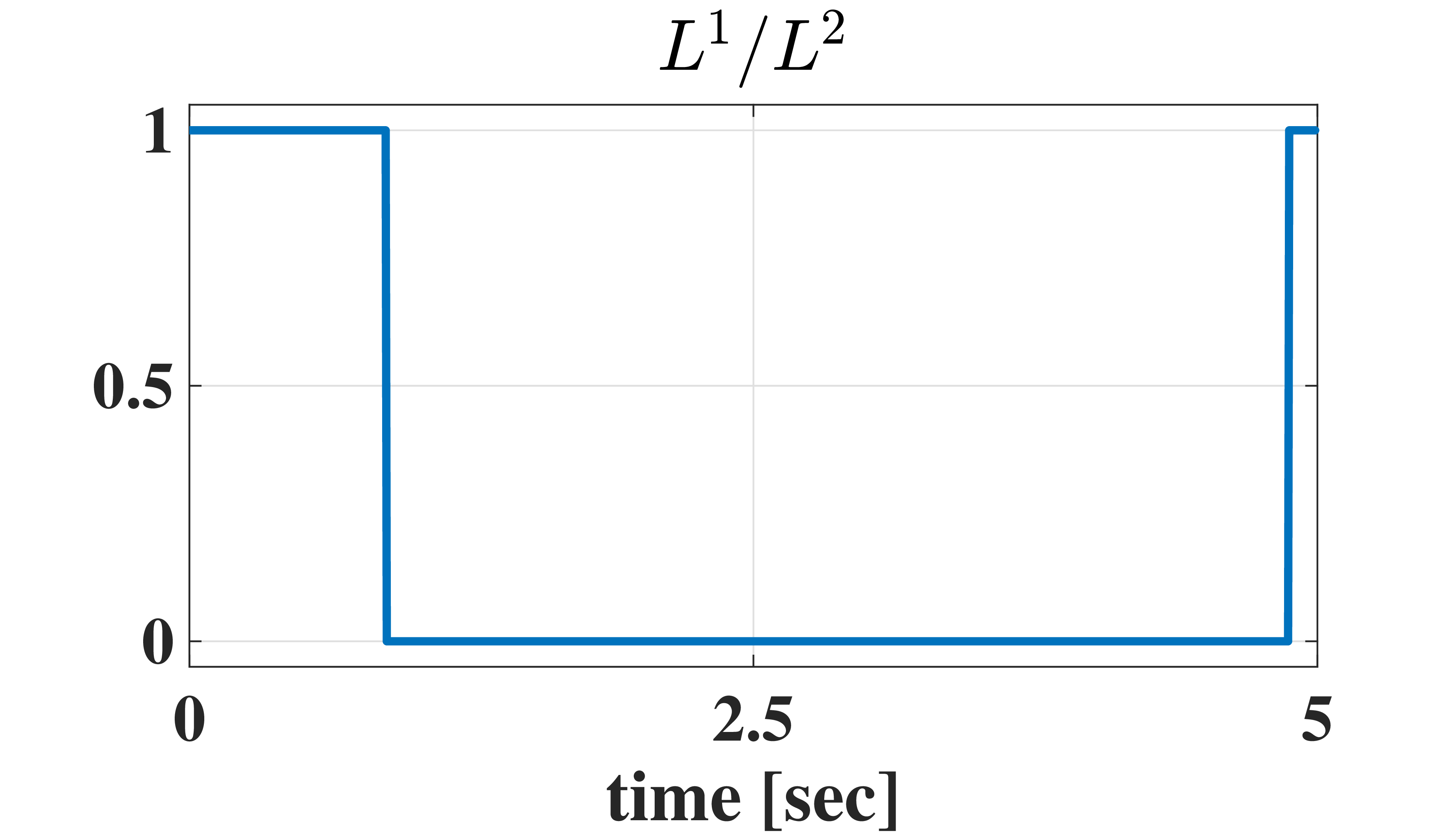}
	\caption{Optimal solutions to various non-convex optimal control problems equivalent to the maximum hands-off control, and the $L^1$ optimal control for comparison.
    The $L^1$ optimal control (top left);
    The $L^p$ optimal control with $(p, \lambda) = (0.5, 0.8)$ (middle left);
    The MCP optimal control with $(\lambda, \alpha) = (1, 0.5)$ (bottom left);
    The SCAD optimal control with $(\lambda, \alpha) = (0.25, 3)$ (top right);
	The LSP optimal control with $(\lambda, \alpha) = (0.1 / \log(1 + 1/\alpha), 10^{-6})$ (middle right);
	The $L^1/L^2$ optimal control with $\lambda = 0.1$ (bottom right).}
\label{fig:ex1_control}
\end{figure} 

\begin{figure}[tb]
  \centering
    \includegraphics[width=0.47\linewidth]{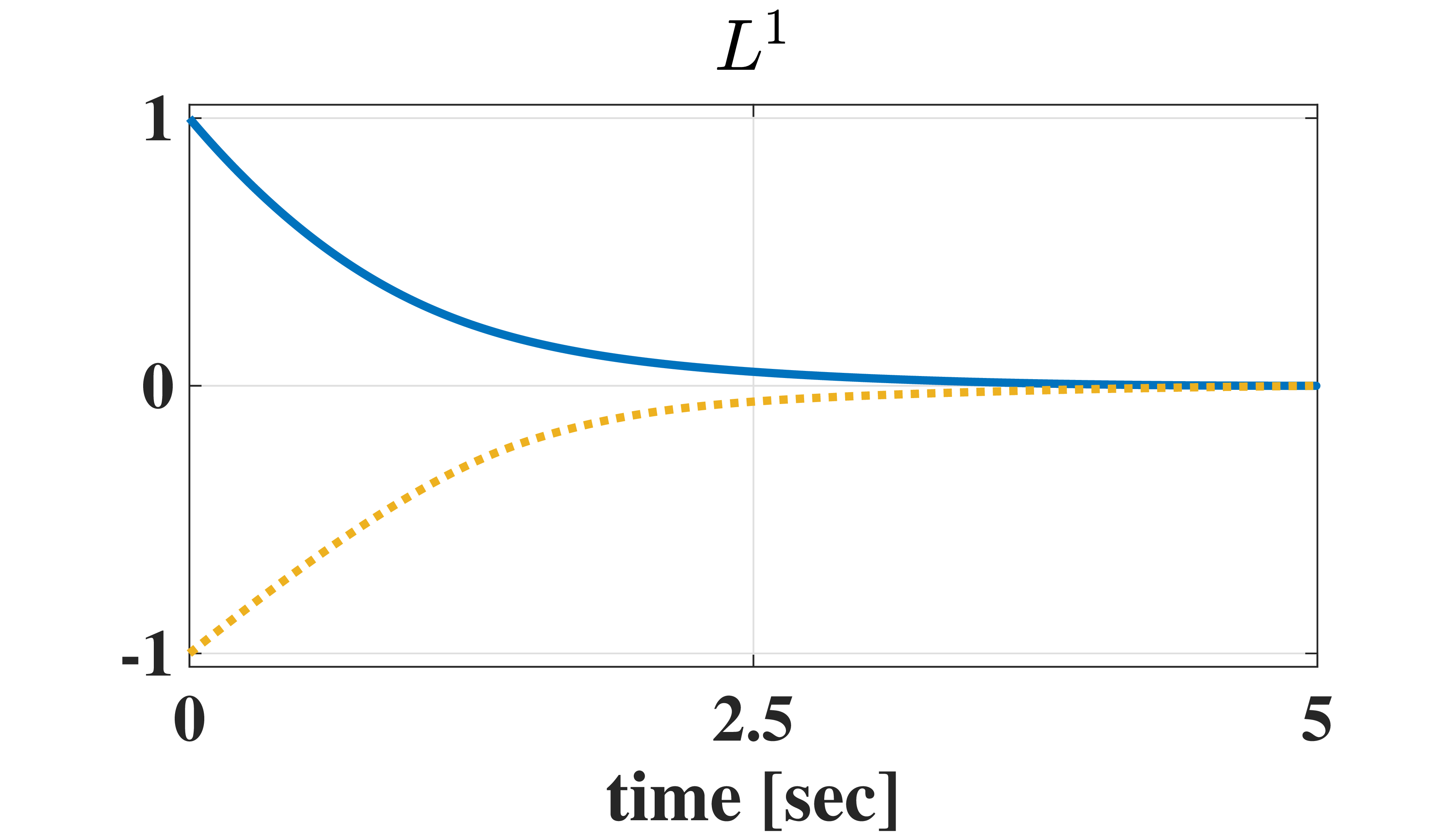}
    \includegraphics[width=0.47\linewidth]{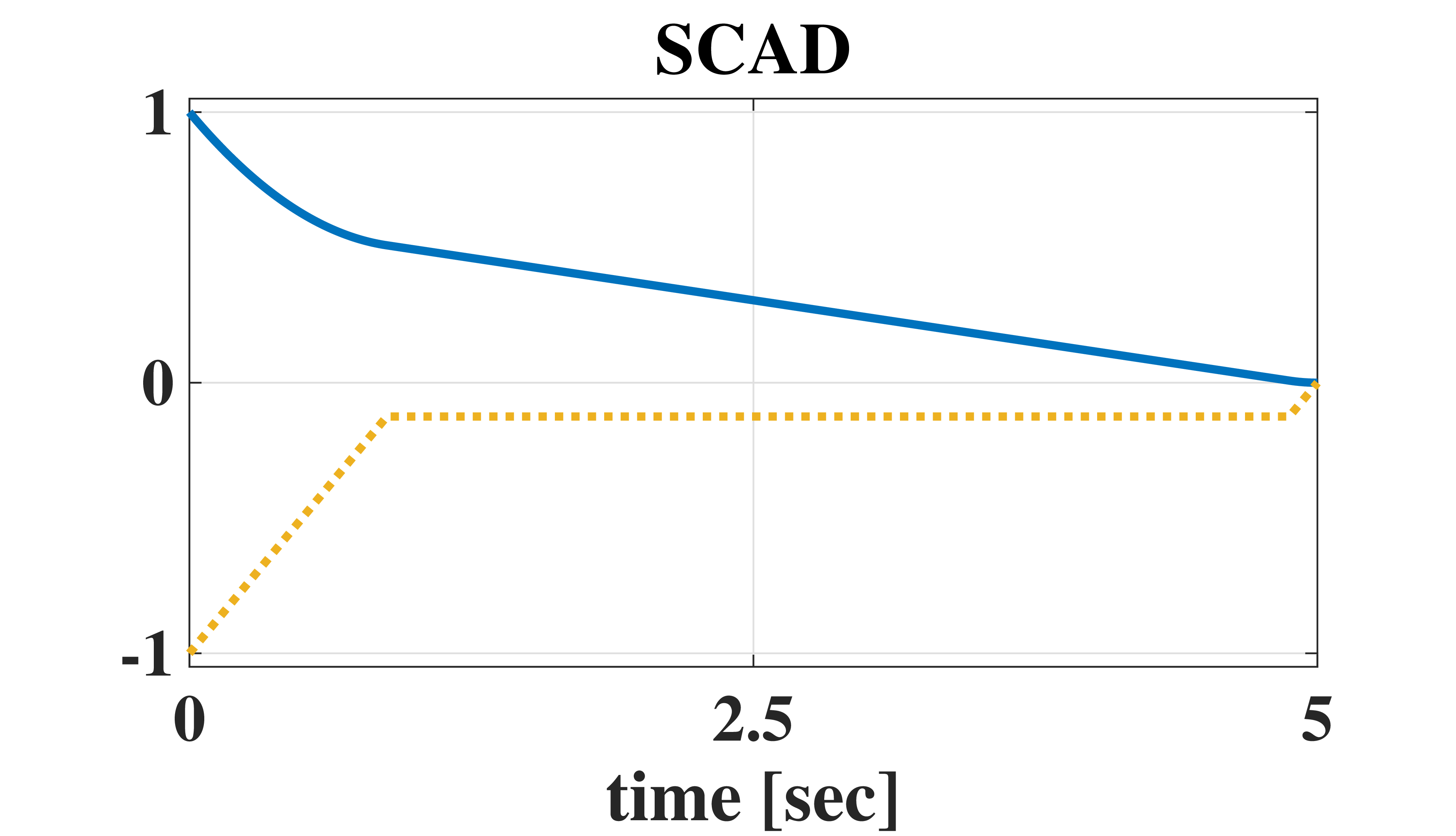}
    \includegraphics[width=0.47\linewidth]{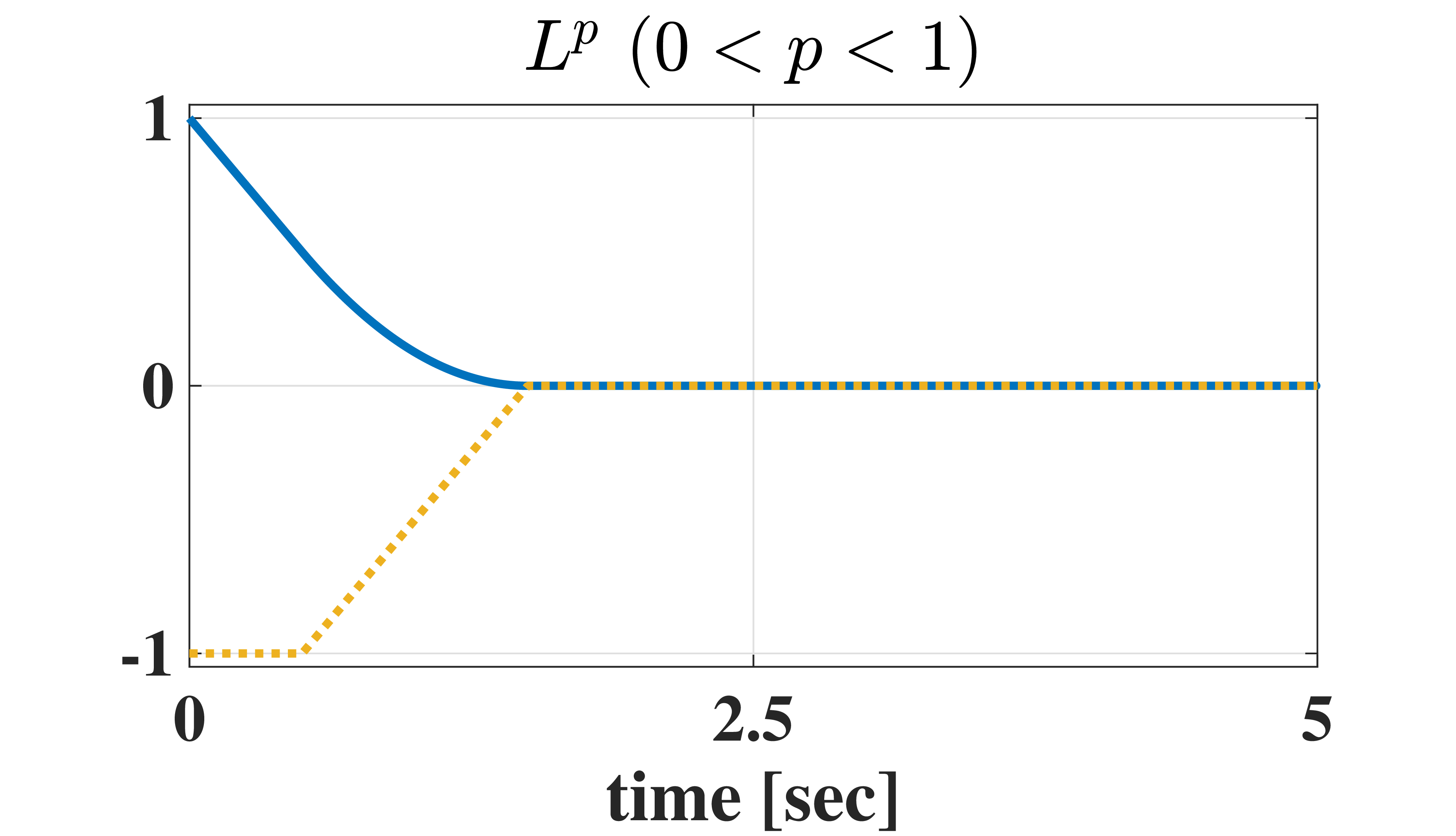}
    \includegraphics[width=0.47\linewidth]{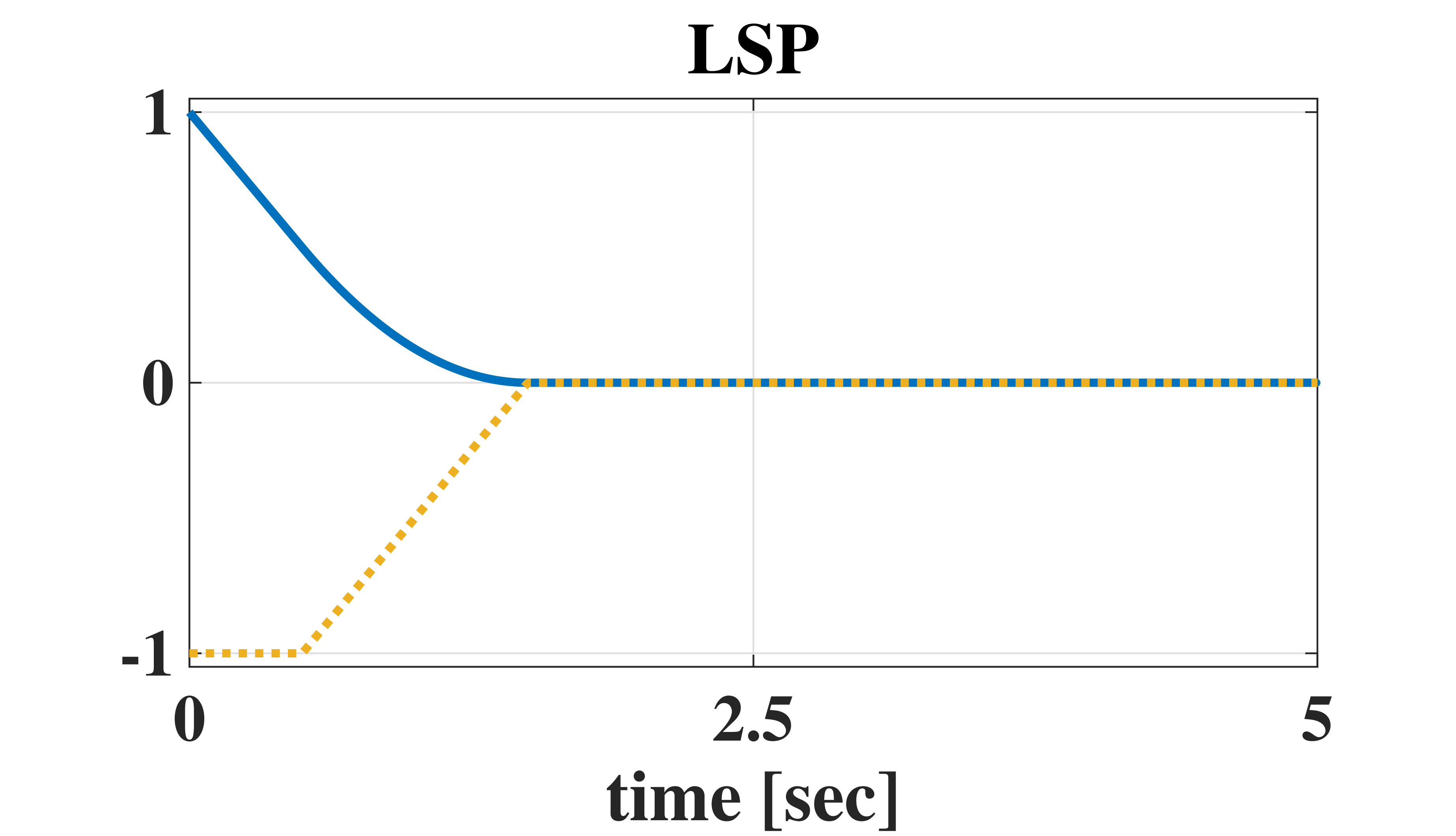}
    \includegraphics[width=0.47\linewidth]{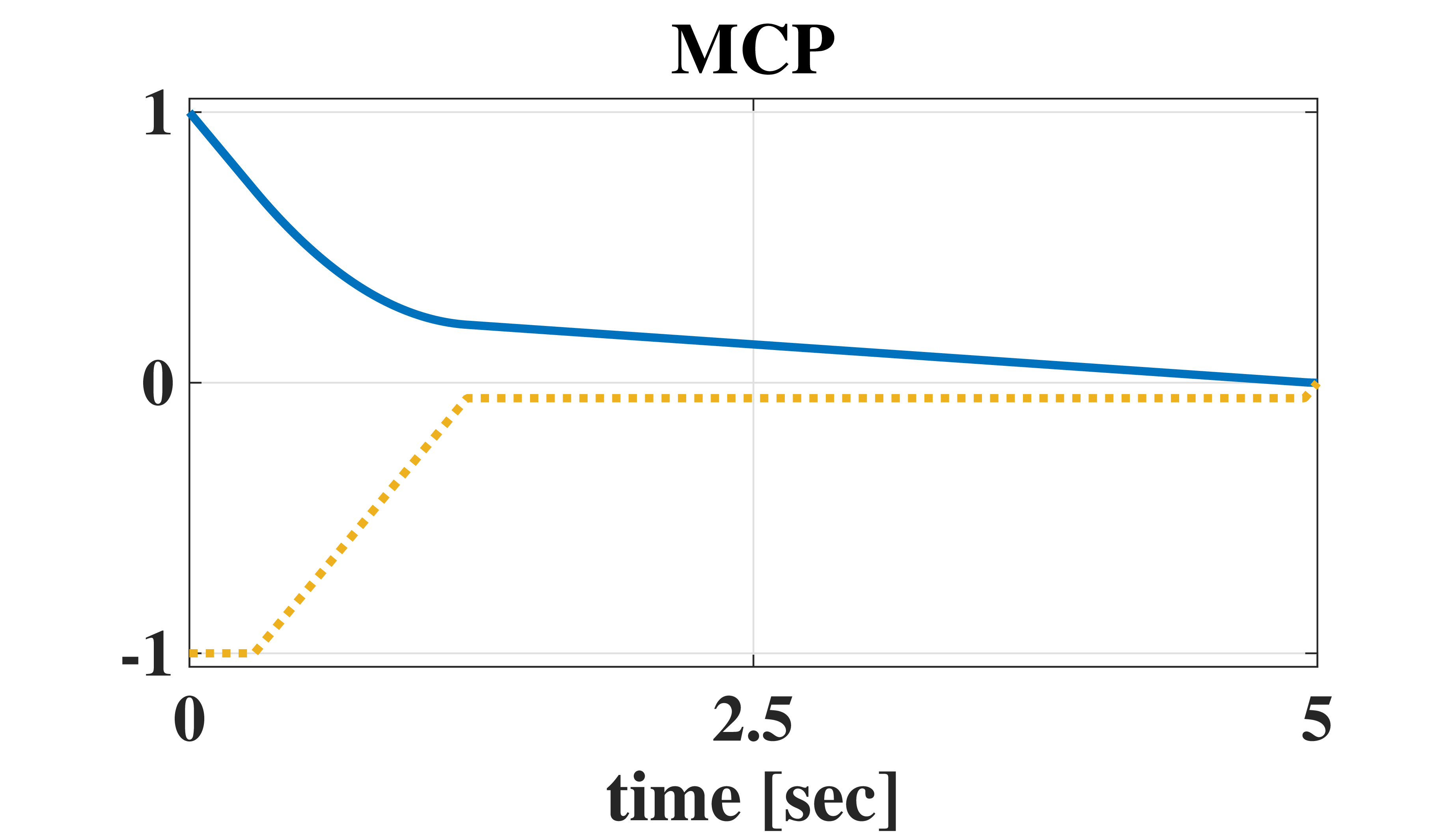}
    \includegraphics[width=0.47\linewidth]{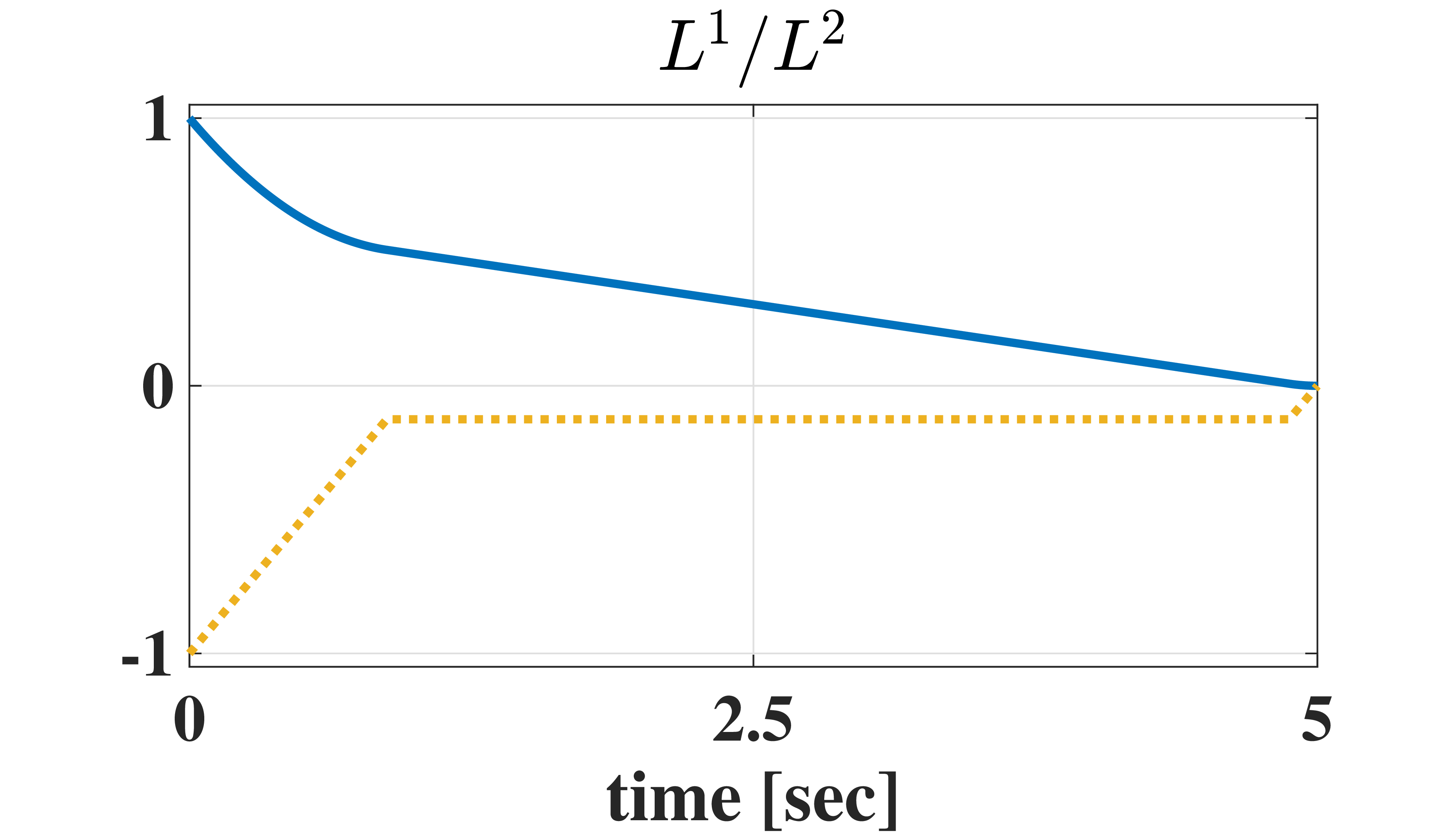}
  \caption{State trajectories, each formed by the optimal solutions in Fig.~\ref{fig:ex1_control}.
    The solid lines and the dotted lines show the state variables $x_1$ and $x_2$ on $[0, T]$, respectively.}
\label{fig:ex1_state}
\end{figure} 

In this section, the results of this paper are confirmed using numerical examples.
The double-integrator system defined by the ordinary differential equation~\eqref{eq:system} is considered, where
\begin{align*}
	A =
	\begin{bmatrix}
		0 & 1\\
		0 & 0
	\end{bmatrix},\quad
	B =
	\begin{bmatrix}
		0\\
		1
	\end{bmatrix}.
\end{align*}
For this system, Problem~\ref{prob:main} with $x_0=\begin{bmatrix}1 &-1\end{bmatrix}^\top$ and $T=5$ is considered.
Here, from~\cite[Theorem 3]{ITKKTAC18}, the set of all maximum hands-off controls is equal to the set of all $L^1$ optimal controls having the bang-off-bang property.
Combining this with~\cite[Control Law 8-3]{AthFal}, the set of all maximum hands-off controls is equal to the set of all inputs $u$ satisfying
\begin{align}\label{eq:ex1_L0}
\begin{split}
	& u(t) \in \{0, 1\} \mbox{~for~almost~all~} t\in[0,T],\\
	& \|u\|_{L^0} = - \xi_{2}, \quad \int_{0}^{T} \int_{0}^{\theta} u(t) dt d\theta = -\xi_1-\xi_2T
\end{split}
\end{align}
where $x_0=\begin{bmatrix} \xi_1 & \xi_2 \end{bmatrix}^\top$.
The first condition is on the bang-off-bang property, 
the second condition on the optimal value, 
and the third condition on  the state constraint.
Therefore, the effectiveness of the non-convex optimal controls can be verified through this example.

The optimal controls introduced in Remark~\ref{rmk:penalties} were computed from Algorithm~\ref{alg:DCA} 
provided that 
$N=1000$,
the $L^1$ optimal control was employed for the initial guess $z[0]$,
and CVX was used for each convex optimization.
Fig.~\ref{fig:ex1_control} shows the obtained control inputs, and the $L^1$ optimal control for comparison.
Fig.~\ref{fig:ex1_state} shows the corresponding state trajectories.
From these figures, all non-convex optimal controls satisfy~\eqref{eq:ex1_L0} and succeed to yield the maximum hands-off control. 
On the other hand, the $L^1$ approximation method fails to yield a sparse solution.
In fact, the $L^1$ optimal solution is not necessarily sparse as shown in~\cite[Control Law 8-3]{AthFal}.
These confirm the effectiveness of the non-convex approximation method.
Finally, the computation time required to find each optimal solution using a standard computer with a 2.7 GHz Intel Core i7 processor is as follows:
0.2333 sec ($L^1$),
1.1058 sec ($L^p$),  
0.8904 sec (MCP), 
0.9106 sec (SCAD),
1.3196 sec (LSP), 
and 1.0144 sec ($L^1/L^2$).
In this example, a convex optimization subproblem was solved three or four times in the DC algorithm. 
Arising from this, the computation time required to solve the non-convex optimization was longer than that to solve the $L^1$ optimization. 
The author plans to work on the improvement of computation algorithm and the selection of the function $\phi$.

\section{Conclusion}
\label{sec:conclusion}
This paper has analyzed the mathematical relationship between a class of some non-convex optimal control problems and the maximum hands-off control problem for continuous-time linear systems. 
The representation capability of the class is compatible with various penalties appeared in the sparse optimization field. 
This paper has proved that the optimal control problems belonging to the class are always equivalent to the maximum hands-off control problem as a main theoretical contribution. 
This property is critically different from the results of the standard approximation method. 
In the numerical computation of the non-convex problems, DC representation was rendered to the maximum hands-off control problem by confining the function $\phi$ to be convex on $[0, 1]^m$, and the computational algorithm quoting the best-known algorithm in the DC optimization field was rendered. 
Then, its effectiveness was confirmed through numerical examples in which the maximum hands-off control can be analytically described.

The algorithm used this time requires computation of a convex subproblem more than one time. 
This could be a drawback in larger systems.
Also, there is arbitrariness in the choice of $\phi$.
Therefore, the improvement of the algorithm and the selection of the specific function $\phi$ suitable for faster computation are future tasks.

\appendix
\section*{Proof of Proposition~\ref{prop:convex_reformulation}}
Let us denote the set of all pairs $(v, w)$ satisfying the constraints in the problem~\eqref{prob:convex_2} by
$\mathcal{U}_3$, and the set of all optimal solutions to the problem~\eqref{prob:convex_2} by $\mathcal{U}_3^\ast$.
As in Theorem~\ref{thm:L0-L1}, the set of all optimal solutions to Problem~\ref{prob:convex} is denoted by $\mathcal{U}_2^\ast$.
From Theorem~\ref{thm:L0-L1}, the set $\mathcal{U}_2^\ast$ is not empty.
Then, let us take any optimal solution $u^\ast\in\mathcal{U}_2^\ast$,
and define $v^\ast$ and $w^\ast$ corresponding to $u^\ast$ by~\eqref{eq:u_v_w_P2}.

Firstly, let us show the statement~(i), i.e., $(v^\ast, w^\ast)\in\mathcal{U}_3^\ast$.
From the definition, 
$u^\ast = v^\ast - w^\ast$,
$(v^\ast, w^\ast)\in\mathcal{U}_3$,
and $v_j^\ast(t) w_j^\ast(t) = 0$ holds on $[0,T]$ for all $j$.
Here, for any $(v, w)\in\mathcal{U}_3$ satisfying $v_j(t)w_j(t)=0$ on $[0, T]$ for all $j$, 
\begin{equation}\label{eq:vw_0}
	J(v) + J(w) = J(u)
\end{equation}
holds, where $u = v - w$, and the above equality follows from
\begin{align*}
	&|u_j(\cdot)| = |v_j(\cdot)-w_j(\cdot)| = v_j(\cdot) + w_j(\cdot),\\
	&\phi_j(u_j(\cdot)) = \phi_j(v_j(\cdot)-w_j(\cdot)) = \phi_j(v_j(\cdot)) + \phi_j(w_j(\cdot)).
\end{align*}
Hence, 
\begin{equation}\label{eq:cost_P2_equiv}
	J(v^\ast) + J(w^\ast) = J(u^\ast).
\end{equation}

Fix any $(v, w)\in\mathcal{U}_3$ and define 
\begin{align}\label{eq:I_positive}
\begin{split}
	&\mathcal{I}_{j}^{(1)} \triangleq \{t \in [0,T]: v_j(t) > w_j(t) > 0\},\\
	&\mathcal{I}_{j}^{(2)} \triangleq \{t \in [0,T]: w_j(t) > v_j(t) > 0\},\\
	&\mathcal{I}_{j}^{(3)} \triangleq \{t \in [0,T]: v_j(t) = w_j(t) > 0\}.
\end{split}
\end{align}
From Lemma~\ref{lem:existence_bang-off}, 
for each $j=1,2,\dots,m$, there exist functions $\rho_j$ and $\sigma_j$ such that 
\begin{align}
	&\rho_j(t)\in\{0, 1\} ~ \forall t\in \mathcal{I}_j^{(1)}, \quad 
		\sigma_j(t)\in\{0, 1\} ~ \forall t\in \mathcal{I}_j^{(2)}, \label{eq:rho_sigma_1}\\
	&\int_{\mathcal{I}_{j}^{(1)}} 
		\begin{bmatrix} e^{-At} b_j \\ 1 \end{bmatrix} (v_j(t) - w_j(t))dt
	= \int_{\mathcal{I}_{j}^{(1)}} 
		\begin{bmatrix} e^{-At} b_j \\ 1 \end{bmatrix} \rho_j(t) dt,\label{eq:rho_sigma_2}\\
	&\int_{\mathcal{I}_{j}^{(2)}} 
		\begin{bmatrix} e^{-At} b_j \\ 1 \end{bmatrix} (v_j(t) - w_j(t))dt
	= -\int_{\mathcal{I}_{j}^{(2)}} 
		\begin{bmatrix} e^{-At} b_j \\ 1 \end{bmatrix} \sigma_j(t) dt.\label{eq:rho_sigma_3}
\end{align}
Then, the functions $v^\circ$ and $w^\circ$ are defined by
\begin{equation}\label{eq:vw_circ}
	(v_j^\circ(t), w_j^\circ(t))=
		\begin{cases}
			(\rho_j(t), 0), 		& \mbox{~if~} t\in\mathcal{I}_{j}^{(1)},\\
			(0, \sigma_j(t)),	& \mbox{~if~} t\in\mathcal{I}_{j}^{(2)},\\
			(0, 0), 				& \mbox{~if~} t\in\mathcal{I}_{j}^{(3)},\\
			(v_j(t), w_j(t)), 	& \mbox{~otherwise},
		\end{cases}
\end{equation}
where $v_j^\circ(t)$ and $w_j^\circ(t)$ are the $j$th component of $v^\circ(t)$ and $w^\circ(t)$, respectively.
Now, 
$(v^\circ, w^\circ)\in\mathcal{U}_3$, and $v^\circ_j(t) w^\circ_j(t) = 0$ holds on $[0, T]$ for all $j$.
Therefore,
\begin{equation}\label{eq:any_vw_P2}
	J(v^\circ) + J(w^\circ) = J(u^\circ)
\end{equation}
holds for $u^\circ = v^\circ - w^\circ$ from~\eqref{eq:vw_0}.
Also,
\begin{align}\label{ineq:I1}
\begin{split}
	&\int_{\mathcal{I}_j^{(1)}}(v_j^\circ(t)+w_j^\circ(t) - \phi_j(v_j^\circ(t))-\phi_j(w_j^\circ(t))) dt \\
	& = \int_{\mathcal{I}_j^{(1)}}\rho_j(t) dt - \phi_j(1) \int_{\mathcal{I}_j^{(1)}} \rho_j(t) dt\\
	& = (1-\phi_j(1)) \int_{\mathcal{I}_j^{(1)}} ( v_j(t)-w_j(t) ) dt\\
	& \leq (1-\phi_j(1)) \int_{\mathcal{I}_j^{(1)}} v_j(t) dt \\
	& \leq \int_{\mathcal{I}_j^{(1)}} ( v_j(t) - \phi_j(v_j(t)) ) dt\\
	& \leq \int_{\mathcal{I}_j^{(1)}} ( v_j(t) - \phi_j(v_j(t)) + w_j(t) - \phi_j(w_j(t))) dt
\end{split}
\end{align}
holds, where \eqref{eq:rho_sigma_1}, \eqref{eq:rho_sigma_2}, $w_j(t)>0$ on $\mathcal{I}_j^{(1)}$, and the assumption~(A3) were used.
In the same way, 
\begin{align}\label{ineq:I2}
\begin{split}
	&\int_{\mathcal{I}_j^{(2)}}(v_j^\circ(t)+w_j^\circ(t) - \phi_j(v_j^\circ(t))-\phi_j(w_j^\circ(t))) dt\\
	&\leq \int_{\mathcal{I}_j^{(2)}} ( v_j(t) - \phi_j(v_j(t)) + w_j(t) - \phi_j(w_j(t))) dt.
\end{split}
\end{align}
Also, from the assumption~(A3),
\begin{align}\label{ineq:I3}
\begin{split}
	&\int_{\mathcal{I}_j^{(3)}}(v_j^\circ(t)+w_j^\circ(t) - \phi_j(v_j^\circ(t))-\phi_j(w_j^\circ(t))) dt\\
	&=0\\
	&\leq \int_{\mathcal{I}_j^{(3)}} ( v_j(t) - \phi_j(v_j(t)) + w_j(t) - \phi_j(w_j(t))) dt.
\end{split}
\end{align}
Therefore, 
\begin{equation}\label{eq:any_vw}
	J(v^\circ) + J(w^\circ) \leq J(v) + J(w).
\end{equation}
Hence,
\begin{align}
	J(v^\ast) + J(w^\ast) 
	&= J(u^\ast)\label{eq:P2_give_P3_1}\\
	&\leq J(u^\circ)\label{eq:P2_give_P3_2}\\
	&= J(v^\circ) + J(w^\circ) \label{eq:P2_give_P3_3}\\
	&\leq J(v) + J(w)\label{eq:P2_give_P3_4}
\end{align}
holds, where
\eqref{eq:P2_give_P3_1} follows from~\eqref{eq:cost_P2_equiv},
\eqref{eq:P2_give_P3_2} from $u^\ast\in\mathcal{U}_2^\ast$ and $u^\circ\in\mathcal{U}$,
\eqref{eq:P2_give_P3_3} from~\eqref{eq:any_vw_P2},
and \eqref{eq:P2_give_P3_4} from~\eqref{eq:any_vw}.
From these and $(v^\ast, w^\ast)\in\mathcal{U}_3$, $(v^\ast, w^\ast)\in\mathcal{U}_3^\ast$ holds.

Next, let us show the statement~(ii).
The set $\mathcal{U}_3^\ast$ is not empty from $(v^\ast, w^\ast)\in\mathcal{U}_3^\ast$.
Then, 
let us take any $(\hat{v}, \hat{w})\in\mathcal{U}_3^\ast$,
define the sets 
$\hat{\mathcal{I}}_{j}^{(1)}$, $\hat{\mathcal{I}}_{j}^{(2)}$, $\hat{\mathcal{I}}_{j}^{(3)}$ for $(\hat{v}, \hat{w})$ as in~\eqref{eq:I_positive},
and construct functions $\hat{v}^\circ$ and $\hat{w}^\circ$ from $\hat{v}$ and $\hat{w}$ as in~\eqref{eq:vw_circ}.
Then, 
$(\hat{v}^\circ, \hat{w}^\circ)\in\mathcal{U}_3$, 
$\hat{u}^\circ\triangleq \hat{v}^\circ - \hat{w}^\circ\in\mathcal{U}$, 
and
\begin{align*}
	J(\hat{v}) + J(\hat{w}) 
	&= J(v^\ast) + J(w^\ast) \\
	&= J(u^\ast)\\
	&\leq J(\hat{u}^\circ)\\
	&= J(\hat{v}^\circ) + J(\hat{w}^\circ) \\
	&\leq J(\hat{v}) + J(\hat{w}),
\end{align*}
which follows from $(v^\ast, w^\ast)\in\mathcal{U}_3^\ast$, \eqref{eq:P2_give_P3_1}, \eqref{eq:P2_give_P3_2}, \eqref{eq:P2_give_P3_3}, \eqref{eq:P2_give_P3_4}.
This gives 
\begin{align}\label{eq:hat_ast_circ}
	J(\hat{v}) + J(\hat{w}) 
	= J(u^\ast)
	= J(\hat{v}^\circ) + J(\hat{w}^\circ).
\end{align}
Then, the inequality holds as the equality in~\eqref{ineq:I1}, \eqref{ineq:I2}, and \eqref{ineq:I3},
where 
$(v_j, w_j, v_j^\circ, w_j^\circ, \mathcal{I}_j^{(l)})$ are replaced by 
$(\hat{v}_j, \hat{w}_j, \hat{v}_j^\circ, \hat{w}_j^\circ, \hat{\mathcal{I}}_j^{(l)})$.
From this, $\mu(\hat{\mathcal{I}}_j^{(1)})=\mu(\hat{\mathcal{I}}_j^{(2)})=\mu(\hat{\mathcal{I}}_j^{(3)})=0$ holds for all $j$.
In other words, $\hat{v}_j(t) \hat{w}_j(t)=0$ on $[0,T]$ for all $j$,
and therefore for $\hat{u}\triangleq \hat{v}-\hat{w}\in\mathcal{U}$, 
$J(\hat{v}) + J(\hat{w}) = J(\hat{u})$
from~\eqref{eq:vw_0}.
From this and~\eqref{eq:hat_ast_circ}, $J(\hat{u}) = J(u^\ast)$ holds,
and therefore $\hat{u} \in \mathcal{U}_2^\ast$.
This completes the proof.

\bibliographystyle{IEEEtran}
\bibliography{IEEEabrv,main}

\begin{thebibliography}{10}
\providecommand{\url}[1]{#1}
\csname url@samestyle\endcsname
\providecommand{\newblock}{\relax}
\providecommand{\bibinfo}[2]{#2}
\providecommand{\BIBentrySTDinterwordspacing}{\spaceskip=0pt\relax}
\providecommand{\BIBentryALTinterwordstretchfactor}{4}
\providecommand{\BIBentryALTinterwordspacing}{\spaceskip=\fontdimen2\font plus
\BIBentryALTinterwordstretchfactor\fontdimen3\font minus
  \fontdimen4\font\relax}
\providecommand{\BIBforeignlanguage}[2]{{%
\expandafter\ifx\csname l@#1\endcsname\relax
\typeout{** WARNING: IEEEtran.bst: No hyphenation pattern has been}%
\typeout{** loaded for the language `#1'. Using the pattern for}%
\typeout{** the default language instead.}%
\else
\language=\csname l@#1\endcsname
\fi
#2}}
\providecommand{\BIBdecl}{\relax}
\BIBdecl

\bibitem{ZhaXuYan15}
Z.~Zhang, Y.~Xu, J.~Yang, X.~Li, and D.~Zhang, ``A survey of sparse
  representation: algorithms and applications,'' \emph{IEEE Access}, vol.~3,
  pp. 490--530, 2015.

\bibitem{Nat95}
B.~K. Natarajan, ``Sparse approximate solutions to linear systems,'' \emph{SIAM
  Journal on Computing}, vol.~24, no.~2, pp. 227--234, 1995.

\bibitem{CanPla09}
E.~J. Cand{\`e}s and Y.~Plan, ``Near-ideal model selection by $\ell_1$
  minimization,'' \emph{Annals of Statistics}, vol.~37, no.~5A, pp. 2145--2177,
  2009.

\bibitem{FanLi01}
J.~Fan and R.~Li, ``Variable selection via nonconcave penalized likelihood and
  its oracle properties,'' \emph{Journal of the American Statistical
  Association}, vol.~96, no. 456, pp. 1348--1360, 2001.

\bibitem{Zha10}
C.-H. Zhang, ``Nearly unbiased variable selection under minimax concave
  penalty,'' \emph{Annals of Statistics}, vol.~38, no.~2, pp. 894--942, 2010.

\bibitem{FraFri93}
L.~E. Frank and J.~H. Friedman, ``A statistical view of some chemometrics
  regression tools,'' \emph{Technometrics}, vol.~35, no.~2, pp. 109--135, 1993.

\bibitem{CanWakBoy08}
E.~J. Cand{\`e}s, M.~B. Wakin, and S.~P. Boyd, ``Enhancing sparsity by
  reweighted $\ell_1$ minimization,'' \emph{Journal of Fourier Analysis and
  Applications}, vol.~14, pp. 877--905, 2008.

\bibitem{Cha07}
R.~Chartrand, ``Exact reconstruction of sparse signals via nonconvex
  minimization,'' \emph{IEEE Signal Processing Letters}, vol.~14, no.~10, pp.
  707--710, 2007.

\bibitem{ChaSta08}
R.~Chartrand and V.~Staneva, ``Restricted isometry properties and nonconvex
  compressive sensing,'' \emph{Inverse Problems}, vol.~24, no.~3, p. 035020,
  2008.

\bibitem{TrzMan08}
J.~Trzasko and A.~Manduca, ``Highly undersampled magnetic resonance image
  reconstruction via homotopic $\ell_0$-minimization,'' \emph{IEEE Transactions
  on Medical Imaging}, vol.~28, no.~1, pp. 106--121, 2008.

\bibitem{TraWeb19}
H.~Tran and C.~Webster, ``A class of null space conditions for sparse recovery
  via nonconvex, non-separable minimizations,'' \emph{Results in Applied
  Mathematics}, vol.~3, p. 100011, 2019.

\bibitem{YinMirPal20}
J.~Ying, J.~V.~M. Cardoso, and D.~Palomar, ``Nonconvex sparse graph learning
  under {L}aplacian constrained graphical model,'' \emph{Advances in Neural
  Information Processing Systems}, vol.~33, pp. 7101--7113, 2020.

\bibitem{WooCha16}
J.~Woodworth and R.~Chartrand, ``Compressed sensing recovery via nonconvex
  shrinkage penalties,'' \emph{Inverse Problems}, vol.~32, no.~7, p. 075004,
  2016.

\bibitem{YinLouHe15}
P.~Yin, Y.~Lou, Q.~He, and J.~Xin, ``Minimization of $\ell_{1-2}$ for
  compressed sensing,'' \emph{SIAM Journal on Scientific Computing}, vol.~37,
  no.~1, pp. A536--A563, 2015.

\bibitem{SouBlaAub15}
E.~Soubies, L.~Blanc-F{\'e}raud, and G.~Aubert, ``A continuous exact $\ell_0$
  penalty ({CEL0}) for least squares regularized problem,'' \emph{SIAM Journal
  on Imaging Sciences}, vol.~8, no.~3, pp. 1607--1639, 2015.

\bibitem{ManKutBru21}
K.~Manohar, J.~N. Kutz, and S.~L. Brunton, ``Optimal sensor and actuator
  selection using balanced model reduction,'' \emph{IEEE Transactions on
  Automatic Control}, vol.~67, no.~4, pp. 2108--2115, 2022.

\bibitem{GomHee15}
T.~Gommans and W.~Heemels, ``Resource-aware {MPC} for constrained nonlinear
  systems: A self-triggered control approach,'' \emph{Systems \& Control
  Letters}, vol.~79, pp. 59--67, 2015.

\bibitem{CarRusHes15}
L.~R.~G. Carrillo, W.~J. Russell, J.~P. Hespanha, and G.~E. Collins, ``State
  estimation of multiagent systems under impulsive noise and disturbances,''
  \emph{IEEE Transactions on Control Systems Technology}, vol.~23, no.~1, pp.
  13--26, 2015.

\bibitem{NagQueNes16}
M.~Nagahara, D.~E. Quevedo, and D.~Ne\v{s}i\'{c}, ``Maximum hands-off control:
  a paradigm of control effort minimization,'' \emph{IEEE Transactions on
  Automatic Control}, vol.~61, no.~3, pp. 735--747, 2016.

\bibitem{Nag2020}
M.~Nagahara, \emph{Sparsity methods for systems and control}.\hskip 1em plus
  0.5em minus 0.4em\relax now publishers, 2020.

\bibitem{ITKKTAC18}
T.~Ikeda and K.~Kashima, ``On sparse optimal control for general linear
  systems,'' \emph{IEEE Transactions on Automatic Control}, vol.~64, no.~5, pp.
  2077--2083, 2018.

\bibitem{ChaNagQueRao16}
D.~Chatterjee, M.~Nagahara, D.~E. Quevedo, and K.~M. Rao, ``Characterization of
  maximum hands-off control,'' \emph{Systems \& Control Letters}, vol.~94, pp.
  31--36, 2016.

\bibitem{NagOstQue16}
M.~Nagahara, J.~{\O}stergaard, and D.~E. Quevedo, ``Discrete-time hands-off
  control by sparse optimization,'' \emph{{EURASIP} Journal on Advances in
  Signal Processing}, vol.~76, pp. 1--8, 2016.

\bibitem{ItoIkeKas21}
K.~Ito, T.~Ikeda, and K.~Kashima, ``Sparse optimal stochastic control,''
  \emph{Automatica}, vol. 125, p. 109438, 2021.

\bibitem{AthFal}
M.~Athans and P.~L. Falb, \emph{Optimal Control}.\hskip 1em plus 0.5em minus
  0.4em\relax Dover Publications, 1966.

\bibitem{Zhan10}
T.~Zhang, ``Analysis of multi-stage convex relaxation for sparse
  regularization.'' \emph{Journal of Machine Learning Research}, vol.~11,
  no.~3, 2010.

\bibitem{HerLas}
H.~Hermes and J.~P. Lasalle, \emph{Function Analysis and Time Optimal
  Control}.\hskip 1em plus 0.5em minus 0.4em\relax Academic Press, 1969.

\bibitem{Ste}
R.~F. Stengel, \emph{Optimal Control and Estimation}.\hskip 1em plus 0.5em
  minus 0.4em\relax Dover Publications, 1994.

\bibitem{LePha18dc}
H.~A.~L. Thi and T.~P. Dinh, ``{DC} programming and {DCA}: thirty years of
  developments,'' \emph{Mathematical Programming}, vol. 169, no.~1, pp. 5--68,
  2018.

\bibitem{TaoAn97}
P.~D. Tao and L.~H. An, ``Convex analysis approach to {DC} programming: theory,
  algorithms and applications,'' \emph{Acta Mathematica Vietnamica}, vol.~22,
  no.~1, pp. 289--355, 1997.

\bibitem{SunYinCheJia18}
T.~Sun, P.~Yin, L.~Cheng, and H.~Jiang, ``Alternating direction method of
  multipliers with difference of convex functions,'' \emph{Advances in
  Computational Mathematics}, vol.~44, pp. 723--744, 2018.

\bibitem{ADMMBoyd}
S.~Boyd, N.~Parikh, E.~Chu, B.~Peleato, and J.~Eckstein, ``Distributed
  optimization and statistical learning via the alternating direction method of
  multipliers,'' \emph{Foundations and Trends in Machine Learning}, vol.~3,
  no.~1, pp. 1--122, 2011.

\bibitem{cvx}
M.~Grant and S.~Boyd, ``{CVX}: Matlab software for disciplined convex
  programming, version 2.1,'' \url{http://cvxr.com/cvx}, Mar. 2014.

\end{thebibliography}

\end{document}